\documentclass[11pt]{amsart}
\usepackage{amsmath, amssymb}
\usepackage{amsfonts}
\usepackage{amsthm}
\usepackage[arrow,matrix,curve,cmtip,ps]{xy}
\usepackage{mathdots}

\allowdisplaybreaks
\newtheorem{theorem}{Theorem}[section]
\newtheorem{lemma}[theorem]{Lemma}

\newtheorem{proposition}[theorem]{Proposition}

\newtheorem*{theorem*}{Theorem}
\newtheorem*{lemmas*}{Lemmas}
\theoremstyle{remark}
\newtheorem{remark}[theorem]{Remark}
\newtheorem{remarks}[theorem]{Remarks}
\theoremstyle{definition}
\newtheorem{definition}[theorem]{Definition}
\newtheorem{notation}[theorem]{Notation}
\newtheorem{example}[theorem]{Example}
\newtheorem{examples}[theorem]{Examples}

\numberwithin{equation}{section}

\newcommand{\N}{\mathbb{N}}

\newcommand{\C}{\mathbb{C}}
\newcommand{\T}{\mathbb{T}}

\newcommand{\aut}{\operatorname{Aut}}

\def\Hom{\operatorname{Mor}}
\def\Obj{\operatorname{Obj}}
\def\id{\operatorname{id}}

\def\pec(#1,#2,#3){[#1;(#2,#3)]}
\def\vcc(#1,#2){[#1;#2]}
\def\shift(#1,#2){\sigma^{#1} {#2}}
\def\Lbar{\overline{\Lambda}}
\def\Lmin{\Lambda^{\min}}
\def\Lbarmin{\overline{\Lambda}^{\min}}
\def\FE{\operatorname{{\mathcal{FE}}}}
\def\Linf{\Lambda^{\leq\infty}}
\def\cst{C^*}
\def\MCE{\operatorname{MCE}}

\begin{document}

\title{Removing sources from higher-rank graphs}
\author{Cynthia Farthing}
\address{Department of Mathematics\\University of Nebraska -- Lincoln\\
Lincoln, NE 68588-0130} \email{cfarthing2@math.unl.edu}
\date{\today}

\begin{abstract}
For a higher-rank graph $\Lambda$ with sources we detail
a construction that creates a row-finite higher-rank graph
$\overline{\Lambda}$ that does not have sources and contains $\Lambda$ as a subgraph.
Furthermore, when $\Lambda$ is row-finite the Cuntz-Krieger algebra of $\Lambda$, $\cst(\Lambda)$
is a full corner of $\cst(\overline{\Lambda})$, the Cuntz-Krieger
algebra of $\overline{\Lambda}$.
\end{abstract}

\maketitle

\section{Introduction}
Higher-rank graphs are generalizations of directed graphs that were introduced by Kumjian and Pask in  \cite{KP1} who were motivated by the $\cst$-algebras of buildings that were studied by Robertson and Steger in \cite{RobSt1,RobSt2,RobSt3}.  In this paper, we extend a higher-rank graph with sources to another higher-rank graph that has no sources.  We will do this in such a way that the $\cst$-algebras of the graphs are strongly Morita equivalent, thereby removing one of the technical difficulties encountered when working with higher-rank graphs.

A higher-rank graph can be viewed as a union of $k$ directed graphs with the same vertex set, where the edges of the different graphs are painted with $k$ different colors.  A higher-rank graph also includes a \emph{factorization property} that dictates how the edges of different colors fit together to form paths.  More precisely, a higher-rank graph $\Lambda$, is a countable category together with a degree functor $d:\Lambda\to\N^k$ which satisfies the factorization property: for every $\lambda\in\Lambda$ and $m,n\in\N^k$ such that $d(\lambda)=m+n$, there are unique elements $\mu,\nu\in\Lambda$ such that $\lambda=\mu\nu$, $d(\mu)=m$ and $d(\nu)=n$. The rank of $\Lambda$ is $k$, and therefore, $\Lambda$ is also called a $k$-graph.  The $\cst$-algebras of higher-rank graphs include the $\cst$-algebras associated to directed graphs which have been the focus of much attention in recent years.  (See \cite{Rae} for a detailed account of graph $\cst$-algebras. We will use the conventions established in \cite{Rae} when discussing directed graphs.)

The development of the $\cst$-algebras of higher-rank graphs has progressed in a manner similar to that of the $\cst$-algebras associated with directed graphs.  The $\cst$-algebras of directed graphs were first defined in terms of groupoids \cite{KPRR}.  Next, in \cite{BPRS}, the graph $\cst$-algebra is realized as the universal $\cst$-algebra generated by a collection of projections and partial isometries satisfying certain relations.  Both of these methods required that the directed graphs be \emph{row-finite}, that is, each vertex has finitely many edges pointing toward it.  The groupoid techniques also required that the directed graph did not have any \emph{sources}. (A source is a vertex that does not have any edges pointing toward it.)  In \cite{FLR}, the $\cst$-algebra of an arbitrary directed graph was defined as a universal $\cst$-algebra.  Using a method similar to that used in \cite{KPRR} for directed graphs, Kumjian and Pask realized the $\cst$-algebra of a higher-rank graph to be the $\cst$-algebra of a groupoid associated to the higher-rank graph.  Therefore, they also required that the higher-rank graphs be {row-finite} and have {no sources} (Definitions \ref{def:rowfinite} and \ref{def:source}).  Raeburn, Sims and Yeend in \cite{RSY1} defined, in a universal way,  the $\cst$-algebras for a class of higher-rank graphs known as \emph{locally convex} $k$-graphs.  Later, they extended their definition to include the $\cst$-algebras of \emph{finitely aligned} $k$-graphs in \cite{RSY2}.  Finitely aligned $k$-graphs allow for vertices to receive infinitely many edges and appear to be the most general class of $k$-graphs to which a $\cst$-algebra can be associated.

One of the main accomplishments of Drinen and Tomforde in \cite{DT} is the development of the method known as \emph{desingularization}.  If $E$ is a directed graph, possibly with sources and possibly not row-finite, a desingularization of $F$ is a row-finite directed graph without sources that is obtained from $E$.  Furthermore, the $\cst$-algebras associated with $E$ and $F$, $\cst(E)$ and $\cst(F)$, respectively, are Morita equivalent.  Therefore, when studying with the $\cst$-algebras associated to directed graphs, it usually suffices to consider directed graphs that are row-finite and have no sources.  The desingularization method, in addition to providing easier proofs for the uniqueness theorems of for graph $\cst$-algebras, also   led to the description of the ideal structure of graph algebras.  (See also \cite{BHRS}.)

The construction detailed in this paper, which ``removes sources" from a higher-rank graph, will have similar effects on the study of higher-rank graph $\cst$-algebras.  First of all, by transforming an arbitrary row-finite higher-rank graph into a locally convex graph, we will be able to use the Cuntz-Krieger relations from \cite[Definition 3.3]{RSY1} which are much simpler than those used to define the algebras of finitely aligned $k$-graphs (Definition~\ref{def:Toeplitz-Cuntz-Kriegerfamily}).  Also, the construction given here may allow for some of the results that exist for the $\cst$-algebras of row-finite higher-rank graphs without sources to be extended to more general higher-rank graph $\cst$-algebras.  For example, in \cite{Ev}, Evans completely describes the $K$-theory of the $\cst$-algebras associated to row-finite $k$-graphs without sources when $k=2$ and obtains some partial results for $k\geq 3$.  Robertson and Sims give necessary and sufficient conditions describing when the $\cst$-algebra corresponding to a row-finite $k$-graph without sources is simple in \cite{RobSi1}.  Since ideal structure and $K$-theory is preserved under Morita equivalence, it is expected that these results will hold in the more general setting.

Our goal is to produce a desingualrization method for higher-rank graphs that is analogous to the process used for directed graphs.  If a vertex $v$ is a source in a directed graph $E$, then the desingularization process ``adds a head to $v$."  This means we attach a graph of the form
$$ \xymatrix{v&v_1\ar[l]_{e_{v_1}}&v_2\ar[l]_{e_{v_2}}&v_3\ar[l]_{e_{v_3}}&\cdots\ar[l]&v_{n-1}\ar[l]
&v_{n}\ar[l]_{e_{v_n}}&\cdots\ar[l]} $$
to $v$.  This method was used by Bates, et. al. in \cite{BPRS} as well as by Drinen and Tomforde in \cite{DT}.

In a directed graph, adding an edge to a vertex automatically creates another directed graph.  Therefore, dealing with sources in a directed graph is a local problem.  However, in a higher-rank graph, adding an edge of some degree to one vertex will require that several edges of different degrees be added to other vertices to ensure that the factorization property still holds.  Hence, adding edges to a vertex in a higher-rank graph is a global issue.  The method we develop here uses the so-called boundary paths of a higher-rank graph to identify the sources and then extends those  boundary paths in the necessary directions.

This paper is designed as follows.  In Section 2, we define the terminology necessary to discuss the $\cst$-algebra of a finitely aligned $k$-graph.  In Section 3, given a row-finite higher-rank graph $\Lambda$, we construct a row-finite higher-rank graph $\overline{\Lambda}$ that is source free.  We show that the $\cst$-algebra of the original $k$-graph sits naturally inside the $\cst$-algebra of the extended $k$-graph as a full corner.  Section 4 includes examples of 2-graphs with sources and how they are extended to graphs without sources using the method in this paper.

\subsection*{Acknowledgements} This research was part of the author's Ph.D. thesis under the direction of Paul Muhly.  The research was conducted during a year long stay at the University of Newcastle, Australia.  The author wishes to thank Paul Muhly as well as Iain Raeburn, Aidan Sims and Trent Yeend for many helpful discussions about this work.

\section{Preliminaries}\label{sec:preliminaries}

\begin{definition}\label{def:kgraph}
Given $k\in\N$, a \emph{$k$-graph} $(\Lambda,d)$ is a  countable category $\Lambda$ together with a functor $d:\Lambda\to\N^k$, called the \emph{degree functor}, which satisfies the \emph{factorization property}: for every $\lambda\in\Hom(\Lambda)$ and $m,n\in\N^k$ with $d(\lambda)=m+n$, there are unique elements $\mu,\nu\in\Hom(\Lambda)$ such that $\lambda=\mu\nu$, $d(\mu)=m$ and $d(\nu)=n$.
\end{definition}

\begin{notation}
\label{N:basic} \hfill{}
\begin{enumerate}
\item For $n\in\N^k$, let $\Lambda^n=\{\lambda\in\Lambda:d(\lambda)=n\}$.  For $E \subseteq \Lambda$ and $\lambda \in \Lambda$, define
\begin{gather*}
\lambda E =\{\lambda\mu : \mu \in E, r(\mu) = s(\lambda)\} \text{ and } \\
E \lambda =\{\mu\lambda : \mu \in E, s(\mu) = r(\lambda)\}.
\end{gather*}
\item We will use $e_1, e_2, \ldots, e_k$ to denote the usual basis for $\N^k$.  For $m,n\in\N^k$, we denote by $m\vee n$ the coordinate-wise maximum and the coordinate-wise minimum by $m\wedge n$.  operations.  Thus $m+n\wedge p=m+(n\wedge p)$ for
$m,n,p\in\N^k$. For $m,n\in\N^k$, $m\vee n$ is the least element in
$\N^k$ that is greater than or equal to both $m$ and $n$, and
$m\wedge n$ is the greatest element in $\N^k$ that is less than or
equal to both $m$ and $n$.

We will use the convention that $\vee$ and $\wedge$ precede addition and subtraction
in the order of For $m,n,p\in\N^k$, it is
straightforward to show that $(m+p)\wedge (n+p)=(m\wedge n)+p$ and
$(m+p)\vee (n+p)=(m\vee n)+p$.
\item Let $\lambda\in\Lambda$ and let $m$ and $n$ satisfy the
inequality $0\leq m\leq n\leq d(\lambda)$. Then the unique
factorization property guarantees that there are unique paths
${\lambda}_{i}$, $i=1,2,3$, such that $d({\lambda}_{1})=m$,
$d({\lambda}_{2})=n-m$,  $d({\lambda}_{3})=d(\lambda)-n$ and
$\lambda = \lambda_1\lambda_2\lambda_3$. We shall write
${\lambda}(0,m)$ for ${\lambda}_{1}$, ${\lambda}(m,n-m)$ for
${\lambda}_{2}$ and ${\lambda}(n,d(\lambda))$ for ${\lambda}_{3}$.
\end{enumerate}
\end{notation}

\begin{examples}
\label{k-graph examples}\hfill{}
\begin{enumerate}
\item Let $E=(E^0, E^1,r,s)$ be a directed graph.
Let $E^*$ denote the category generated freely over all finite
paths. Let $l:E^*\to \N$ give the length of a path. Then $(E^*,l)$
is a 1-graph.
\item For $m\in(\N\cup\{\infty\})^{k}$, define $\Omega_{k,m}$ to be the
$k$-graph with \[
\Obj(\Omega_{k,m})=\{ p\in\N^{k}:p\leq m\},\]
 \[
\Hom(\Omega_{k,m})=\{(p,q)\in\Obj(\Omega_{k,m})\times\Obj(\Omega_{k,m}):p\leq q\},\]
 \[
r(p,q)=p,\quad s(p,q)=q,\quad d(p,q)=q-p.\]  Drawn below are
$\Omega_{2,(\infty,\infty)}$ and $\Omega_{2,(1,2)}$. In the diagrams,
edges of degree $(1,0)$ are solid; edges of degree $(0,1)$ are dashed.  In
each diagram $\lambda=((0,2),(1,2))$ while $\mu~=((0,0),(0,1))$.\\

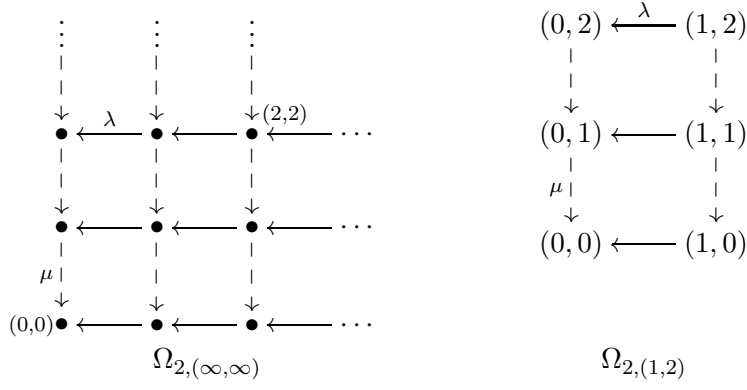
\begin{center}
\begin{tabular}{ccc}
{$\xymatrix
{\vdots\ar@{-->}[d]&\vdots\ar@{-->}[d]&\vdots\ar@{-->}[d]\\
\bullet\ar@{-->}[d]&\bullet\ar@{-->}[d]\ar[l]_{\lambda}&\bullet\ar@{-->}[d]\ar[l]&\cdots\ar[l]_>>>{(2,2)}\\
\bullet\ar@{-->}[d]_{\mu}&\bullet\ar@{-->}[d]\ar[l]&\bullet\ar@{-->}[d]\ar[l]&\cdots\ar[l]\\
\llap{$\scriptstyle (0,0)$}\bullet&\bullet\ar[l]&\bullet\ar[l]&\cdots\ar[l]}
$}&\hspace{.5in}\ &
{$
\xymatrix{(0,2)\ar@{-->}[d]&(1,2)\ar[l]_{\lambda}\ar@{-->}[d]\\
(0,1)\ar@{-->}[d]_{\mu}&(1,1)\ar[l]\ar@{-->}[d]\\
(0,0)&(1,0)\ar[l]}$}\\
$\Omega_{2,(\infty,\infty)}$&&$\Omega_{2,(1,2)}$
\end{tabular}
\end{center}
\end{enumerate}
\end{examples}

\begin{definition}\label{def:rowfinite}
A $k$-graph $(\Lambda, d)$ is \emph{row-finite} if $v\Lambda^n$ is at most finite for all $v\in\Lambda^0$ and $n\in\N^k$.
\end{definition}

\begin{definition}\label{def:source}
A vertex $v\in\Lambda^0$ is a \emph{source} if $v\Lambda^n=\emptyset$ for some $n\in\N^k$.
\end{definition}

\begin{definition}\label{def:minimalcommonextension}
For $\lambda,\mu\in\Lambda$,  if there exist
$\alpha,\beta\in\Lambda$ such that $\lambda\alpha=\mu\beta$ and
$d(\lambda\alpha)=d(\lambda)\vee d(\mu)$, then $\lambda\alpha$ is
called a \emph{minimal common extension of $\lambda$ and $\mu$}.
Define
$$\Lmin(\lambda,\mu)=\{(\alpha,\beta)\in \Lambda\times\Lambda:
\lambda\alpha=\mu\beta \text{ and } d(\lambda\alpha)=d(\lambda)\vee d(\mu)\}.$$
\end{definition}

\begin{definition}\label{def:finitelyaligned}
A $k$-graph $(\Lambda,d)$ is \emph{finitely aligned} if
$\Lmin(\lambda,\mu)$ is at most finite for all
$\lambda,\mu\in\Lambda$.
\end{definition}

\begin{remark}\label{rmk:differences between 1 and k graphs}
Definitions \ref{def:source} and \ref{def:finitelyaligned}
highlight some key differences between 1-graphs and $k$-graphs for
$k\geq 2$.  First of all, in the directed
graph setting, a source is a vertex $v$ for which
$v\Lambda^1=\emptyset$, or equivalently, if $v\Lambda=\{v\}$.  However, a vertex in a 1-graph is a source in the sense of Definition \ref{def:source} if  there exists a path $\lambda\in v\Lambda$ such that
$s(\lambda)\Lambda=\{s(\lambda)\}$. This is not the case for
arbitrary $k$-graphs.  Consider the graph $\Omega_{2,(\infty,1)}$
drawn here.

$$
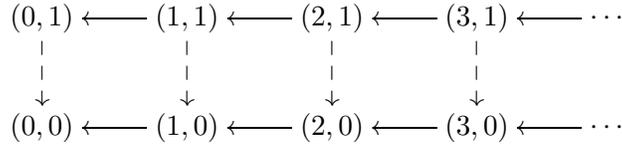
\xymatrix{
(0,1)\ar@{-->}[d]&(1,1)\ar[l]\ar@{-->}[d]&(2,1)\ar[l]\ar@{-->}[d]&(3,1)\ar[l]\ar@{-->}[d]&\cdots\ar[l]\\
(0,0)&(1,0)\ar[l]&(2,0)\ar[l]&(3,0)\ar[l]&\cdots\ar[l]}
$$

\noindent Each of the vertices $(m,1)$, $m\in\N$ is a source  since
$(m,1)\Omega_{2,(\infty,1)}^{e_2}=\emptyset$.  However, there is no
vertex $v\in\Omega_{2,(\infty,1)}^0$ with
$v\Omega_{2,(\infty,1)}=\{v\}$.  The difference is that in a $k$-graph for $k\geq 2$, vertices can
be sources in some directions, but not in all.  Secondly, if
$\Lambda$ is a 1-graph and $\lambda,\mu\in\Lambda$, the only way
two paths can have a minimal common extension is if one
path is a subpath of the other. Therefore, the set
$\Lmin(\lambda,\mu)$ is either empty or a singleton. Consequently, any 1-graph is
finitely aligned.
\end{remark}

\begin{definition}\label{def:exhaustive}
Let $(\Lambda, d)$ be a $k$-graph; let $v\in\Lambda^0$ and $E\subset
v\Lambda$.  We say that $E$ is \emph{exhaustive} if for every
$\mu\in v\Lambda$ there exists a $\lambda\in E$ such that
$\Lmin(\lambda,\mu)\neq\emptyset$. We denote the set of all
\emph{finite exhaustive subsets} of $\Lambda$ by $\FE(\Lambda)$.
\end{definition}

\begin{examples}\label{ex:exhaustive}
\hfill{}
\begin{enumerate}
\item For all $m\in(\N\cup\{\infty\})^k$ and $v\in\Omega_{k,m}^0$,
any nonempty finite subset of $v\Omega_{k,m}$ is a finite exhaustive set.
\item Consider the $k$-graph $\Lambda$ below:

$$
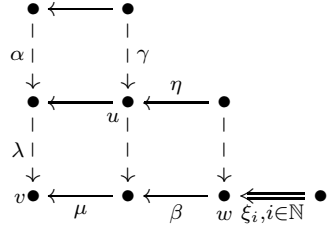
\xymatrix{\bullet\ar@{-->}[d]_{\alpha}  & \bullet\ar[l]\ar@{-->}[d]^{\gamma}\\
\bullet \ar@{-->}[d]_{\lambda}&\bullet \ar[l]\ar@{-->}[d]_<{u}&\bullet \ar[l]_{\eta}\ar@{-->}[d]\\
\llap{$\scriptstyle v$}\bullet &\bullet \ar[l]^{\mu}&\smash{\underset{w}{\bullet}}\vphantom{\bullet} \ar[l]^{\beta}&\bullet\ar@{=>}[l]^{\xi_i, i\in\N}}$$

\noindent Dashed edges represent edges of degree $(0,1)$ and solid
edges represent edges of degree $(1,0)$.  The edges $\xi_i$ where
$i\in\N$ each have degree $(1,0)$. Any finite exhaustive subset of
$w\Lambda$ must contain $w$.  The set $\{\mu\}$ is a finite
exhaustive subset of $v\Lambda$, whereas $\{\lambda\}$ is not
because $\Lmin(\lambda,\mu\beta\xi_i)=\emptyset$ for any $i\in\N$.
\end{enumerate}
\end{examples}

\begin{definition}\label{def:Toeplitz-Cuntz-Kriegerfamily}
Let $(\Lambda, d)$ be a finitely aligned $k$-graph.  A
\emph{Toeplitz-Cuntz-Krieger  $\Lambda$-family} in a $C^*$-algebra
$B$ consists of a family of partial isometries
$\{t_\lambda:\lambda\in\Lambda\}$ satisfying the
\emph{Toeplitz-Cuntz-Krieger relations}:
\begin{list}{}{\setlength{\itemindent}{.1in}\setlength{\labelwidth}{6em}}
\item[(TCK1)] $\{t_v:v\in\Lambda^0\}$ is a family of mutually orthogonal projections
\item[(TCK2)] $t_{\lambda\mu}=t_\lambda t_\mu$ for all $\lambda,\mu\in\Lambda$ with $s(\lambda)=r(\mu)$
\item[(TCK3)] $t^*_\lambda t_\mu = \sum_{(\alpha,\beta)\in\Lmin(\lambda,\mu)}t_\alpha t^*_\beta$
for all $\lambda,\mu\in\Lambda$
\end{list}
A \emph{Cuntz-Krieger $\Lambda$-family} in a $C^*$-algebra $B$ is  a
Toeplitz-Cuntz-Krieger $\Lambda$-family that also satisfies

\noindent\hspace{2em}(CK) $\prod_{\lambda\in
E}(t_v-t_{\lambda}t_{\lambda}^*)=0$  for all $v\in\Lambda^0$ and
$E\in v\FE(\Lambda)$.
\end{definition}

Of course, the hypothesis that $({\Lambda},d)$ is finitely aligned
guarantees that the sums in Definition~\ref{def:Toeplitz-Cuntz-Kriegerfamily} are finite
sums, and hence make sense in any $C^{*}$-algebra.

\begin{definition}\label{def:cst Lambda}
Let $(\Lambda,d)$ be a finitely aligned $k$-graph.  The \emph{$\cst$-algebra
of $\Lambda$},
denoted $\cst(\Lambda)$, is the $\cst$-algebra generated by a
universal Cuntz-Krieger $\Lambda$-family
$\{s_\lambda:\lambda\in\Lambda\}$ which is universal if the sense
that if $\{t_\lambda:\lambda\in\Lambda\}$ is a Cuntz-Krieger
$\Lambda$-family in a $\cst$-algebra $B$, then there exists a
$\cst$-homomorphism $\pi:\cst(\Lambda)\to B$ such that
$\pi(s_\lambda)=t_\lambda$ for all $\lambda\in\Lambda$.
\end{definition}
We also call $\cst(\Lambda)$ the Cuntz-Krieger algebra of $\Lambda$.
\begin{definition}\label{def:graph morphism}
Let $(\Lambda_1,d_1)$ and $(\Lambda_2, d_2)$ be $k$-graphs.  A
\emph{graph morphism} is a functor $F:\Lambda_1\to \Lambda_2$ such
that $d_2(F(\lambda))=d_1(\lambda)$ for all $\lambda\in\Lambda$.
\end{definition}

\begin{definition}\label{def:X Lambda}
Let $(\Lambda,d)$ be a $k$-graph.  Define the \emph{path space of
$\Lambda$}  to be the set
$$X_\Lambda=\{x:\Omega_{k,m}\to\Lambda: m\in(\N\cup\{\infty\})^k\text{ and $x$ is a graph morphism}\}.$$
We extend the range and degree map  of $\Lambda$ to $X_\Lambda$ by
defining, for $x:\Omega_{k,m}\to\Lambda$, $r(x)=x(0)$ and $d(x)=m$.
\end{definition}

\begin{remarks}
\hfill{}
\begin{enumerate}
\item The factorization property of $k$-graphs implies that each $x\in
X_\Lambda$ is completely determined by $\{x(0,p): p\leq m\}$: if
$l\leq n\leq p$ and $x(0,p)=\lambda_p$, then
$x(l,n)=\lambda_p(l,n)$.  If $m_i<\infty$ for all
$i\in\{1,2,\ldots,k\}$, then $x(0,m)$ completely determines $x$.

\item The map $\lambda\mapsto x_\lambda$ from $\Lambda$ to $X_\Lambda$
where $x_\lambda$ is the path discussed above,  embeds $\Lambda$ into
$X_\Lambda$.
\end{enumerate}
\end{remarks}

\begin{notation}
Let $x:\Omega_{k,m}\to\Lambda$ be a graph morphism.
\begin{enumerate}
\item For $p\leq m$, define $\sigma^p x:\Omega_{k,m-p}\to\Lambda$ by $\sigma^p x(a,b)=x(a+p,b+p)$ for $a,b\in\N^k$ such that $a\leq b\leq m-p$.
\item For $\lambda\in\Lambda$ such that $s(\lambda)=x(0)$ define $\lambda x:\Omega_{k,m+d(\lambda)}\to\Lambda$ by $(\lambda x)(0,d(\lambda))=\lambda$ and $(\lambda x)(0,p)=\lambda x(0,p-d(\lambda))$ for $p\in\N^k$ such that $d(\lambda)\leq p\leq d(x)+d(\lambda)$.

\end{enumerate}
\end{notation}

\begin{definition}\label{def:locallyconvex}
A $k$-graph $(\Lambda,d)$ is \emph{locally convex}  if whenever
$\lambda\in v\Lambda^{e_i}$ and $\mu\in v\Lambda^{e_j}$ for some
$v\in\Lambda^0$ and $i,j\in\{1,2,\ldots,k\}$ with $i\neq j$, there
exists $\xi\in s(\lambda)\Lambda^{e_j}$ and $\eta\in
s(\mu)\Lambda^{e_i}$.
\end{definition}

\begin{definition}\label{def:Lambda less than paths}
Let $(\Lambda,d)$ be a $k$-graph.  For $q\in\N^k$, define
$$
\Lambda^{\leq q}=\{\lambda\in\Lambda: d(\lambda)\leq q, \text{ and }
s(\lambda)\Lambda^{e_i}=\emptyset \text{ when } d(\lambda)+e_i\leq
q\}.$$
\end{definition}

\begin{examples}\label{ex:locallyconvex}
\hfill{}
\begin{enumerate}
\item For any $m\in (\N\cup\{\infty\})^k$, $\Omega_{k,m}$ is locally convex.  More generally, if $\Lambda$ has no sources, then $\Lambda$ is locally convex since $v\Lambda^{e_i}\neq\emptyset$ for all $v\in\Lambda^0$ and $i\in\{1,2,\ldots,k\}$.
\item The $2$-graph in Example~\ref{ex:exhaustive}~(2) is not locally convex.
For the vertex $u$, we have $\eta\in u\Lambda^{e_1}$ and $\gamma\in u\Lambda^{e_2}$.  However, $s(\eta)\Lambda^{e_2}$ and $s(\gamma)\Lambda^{e_1}$ are both empty.
\end{enumerate}
\end{examples}

\begin{remark}\label{rmk:locallyconvexCKrelation}
Condition (CK) of Definition \ref{def:Toeplitz-Cuntz-Kriegerfamily}
replaced earlier Cuntz-Krieger conditions used for row-finite
$k$-graphs with no sources \cite{KP1} and for locally convex
$k$-graphs \cite{RSY1}. The condition from \cite{RSY1} is
\begin{equation}\tag{CK$^\prime$}
t_v=\sum_{\lambda\in\Lambda^{\leq m}} t_\lambda t_\lambda^* \text{
for all } v\in\Lambda^0 \text{ and } m\in\N^k.
\end{equation}
It is shown in \cite[Appendix B]{RSY2} that the conditions in
Definition \ref{def:Toeplitz-Cuntz-Kriegerfamily} are equivalent to
those in \cite{RSY1} when the $k$-graph is locally convex.
\end{remark}

\begin{definition}\label{def:boundarypath}
Let $(\Lambda,d)$ be a $k$-graph; let $x:\Omega_{k,m}\to\Lambda$ be
a graph morphism in $X_\Lambda$ for some $m\in(\N\cup\{\infty\})^k$.
Then $x$ is a \emph{boundary path} if there exists $n_x\in\N^k$ such
that $n_x\leq m$ and for $p\in\N^k$
$$ (n_x\leq p\leq m,  \text{ and } p_i=m_i)\Rightarrow x(p)\Lambda^{e_i}=\emptyset.$$
We write $\Lambda^{\leq\infty}$ for the collection of all boundary
paths of $\Lambda$.
\end{definition}

Boundary paths are essential to the construction detailed in the next section. We will use the following results about boundary paths.\\
\vspace{12pt}

\begin{lemmas*}\label{lem:boundarypaths}
\hfill{}
\begin{enumerate}
\item \cite[Lemma 2.10]{RSY2} If $x\in\Linf$, then $\shift(x,p)$ and
$\lambda x$ are elements of $\Linf$ for any $p\leq d(x)$ and
$\lambda\in x(0)\Lambda$.
\item \cite[Lemma 2.11]{RSY2} For any
$v\in\Lambda^0$, the set $v\Linf$ is nonempty.
\end{enumerate}
\end{lemmas*}

\section{Removing Sources}\label{sec:removing sources}
In this section,
we will develop a method that extends a finitely aligned $k$-graph with
sources, named $\Lambda$, to a row-finite $k$-graph without sources,
$\overline{\Lambda}$.  When $\Lambda$ is row-finite, $\cst(\Lambda)$ is Morita equivalent
to $\cst(\overline{\Lambda})$. The following theorem is the goal of
this section.

\begin{theorem}\label{thm:nosources}
Let $(\Lambda,d)$ be a row-finite $k$-graph.   Then there exists a
row-finite $k$-graph $(\overline{\Lambda},\overline{d})$ without sources and an
isomorphism $\iota$ of $\Lambda$ onto a subgraph of
$\overline{\Lambda}$ such that the $\cst$-subalgebra of
$C^*(\overline{\Lambda})$ generated by
$\{s_\lambda:\lambda\in\iota{\Lambda}\}$ is a full corner of
$C^*(\overline{\Lambda})$ and is canonically isomorphic to
$C^*(\Lambda)$.
\end{theorem}

We will spend the rest of the section constructing
$\overline{\Lambda}$ and proving Theorem \ref{thm:nosources}.  We
begin by defining two equivalence relations $\sim$ and $\approx$.
The equivalence classes given by $\sim$ correspond to the paths that
will be added to $\Lambda$, and the equivalence classes of $\approx$
correspond to the new vertices.

\begin{definition}\label{def:new vertices}
Let  $V_\Lambda$=\{$(x;m)$\ : $x\in \Lambda^{\leq \infty}$ and $m
\not\leq d(x)$\}.
\end{definition}

The set $V_\Lambda$ extends each element of $\Linf$ in the proper directions.  Notice that the set $V_\Lambda$ is disjoint from $\Lambda^0$ because every vertex  in $\Lambda$ can be written as $x(m)$ for some
$x\in\Linf$ and $m\leq d(x)$.   However, extending each boundary path separately adds many more vertices to $\Lambda$ than
necessary because boundary paths can overlap. An example of such overlap would occur for paths $x,y\in\Linf$ such that $y=\shift(p,x)$ for some $p\leq d(x)$.  To take possible overlap into account, we define the following relation on $V_\Lambda$.

\begin{definition}\label{def:vertexequivalence}
Define a relation $\approx$ on $V_\Lambda$ by: $(x;m)\approx(y;p)$
if
\begin{list}{}{\setlength{\itemindent}{.6em}\setlength{\labelwidth}{2em}}
    \item[(V1)] \label{veqr1}$x(m\wedge d(x))=y(p\wedge d(y))$
    \item[(V2)] \label{veqr2}$m-m\wedge d(x)=p-p\wedge d(y)$
\end{list}
\end{definition}

Condition (V1) ensures that two new vertices are related if they
project down onto the same vertex in $\Lambda$.  Condition (V2)
relates two vertices in $V_\Lambda$ if they are the same ``distance"
from $\Lambda$.

The proof of the next proposition is clear.

\begin{proposition}\label{prop:vertexeqr}
The relation $\approx$ on $V_\Lambda$ is an equivalence relation.
\end{proposition}

\begin{definition}\label{def:new paths}
Let $P_\Lambda=\{(x;(m,n)):x \in \Lambda^{\le \infty}, n \not\le
d(x), \text{ and } m\le n\}$.
\end{definition}

Recall the definition  of $\Omega_{k,m}$ in Example
\ref{k-graph examples} (ii) where paths were denoted by pairs of vertices.
Definition \ref{def:new paths} uses an analogous way to describe the
paths that extend the original $k$-graph.  Since in
Definition~\ref{def:new paths} $n\not\leq d(x)$ but $m$ may or may
not be less than or equal to $d(x)$, we are requiring that the
additional paths start (have source) outside of the original
$k$-graph but may or may not end (have range) in the original
$k$-graph.  Again, the elements of $P_\Lambda$ are paths extending
each boundary path, and therefore, the overlapping of boundary paths
must be taken into account.

\begin{definition}\label{def:pathequivalence}
Define a relation $\sim$ on $P_\Lambda$ by $(x;(m,n))\sim(y;(p.q))$
if
\begin{list}{}{\setlength{\itemindent}{.6em}\setlength{\labelwidth}{2em}}
    \item[(P1)] \label{peqr1}$x(m\wedge d(x),n\wedge d(x))=y(p\wedge d(y),q\wedge d(y))$
    \item[(P2)] \label{peqr2}$m-m\wedge d(x)=p-p\wedge d(y)$
    \item[(P3)] \label{peqr3}$n-m=q-p$
\end{list}
\end{definition}

\begin{proposition}\label{prop:patheqr}
The relation $\sim$ on $P_\Lambda$ is an equivalence relation.
\end{proposition}

Let $\widetilde{P_\Lambda}={P_\Lambda}/\sim$  and
$\widetilde{V_\Lambda}=V_\Lambda/\approx$.  The equivalence classes
of $\widetilde{P_\Lambda}$ will be denoted $[x;(m,n)]$, and the
equivalence classes of $\widetilde{V_\Lambda}$ will be denoted
$[x;m]$.

As mentioned earlier, our  goal is to define a new category
$\overline{\Lambda}$ that extends $\Lambda$.  The elements of
$\widetilde{V_\Lambda}$ will become the additional objects joined to $\Lambda$, and the new
morphisms will be the elements of $\widetilde{P_\Lambda}$.  We now proceed by defining the range and source maps as well as the compostion ($\circ$) and identity ($\id$) functions on
$\widetilde{P_\Lambda}$ that will be used to define the new category.

\begin{definition}\label{def:rangeandsourceextension}
Define  $\widetilde{r}:\widetilde{P_\Lambda} \to
(\widetilde{V_\Lambda}\cup\Lambda^0)$ and
$\widetilde{s}:\widetilde{P_\Lambda} \to \widetilde{V_\Lambda}$ as
follows:

\begin{align*}
& \widetilde{r}([x;(m,n)])=\begin{cases}
x(m) &\text{if $m\leq d(x)$,}\\
\vcc(x,m) &\text{if $m\not\leq d(x)$,}
\end {cases} \\
&\widetilde{s}([x;(m,n)])=[x;n].
\end{align*}
\end{definition}

Notice that  if $(x;(m,n)),(y;(p,q))\in P_\Lambda$, $m\not\leq
d(x)$, $p\not\leq d(y)$ and $(x;(m,n))\sim (y;(p,q))$, then
Condition (P1) implies $x(m\wedge d(x))=y(p\wedge d(y))$.  This
together with Condition (P2) shows that $(x;m)\approx(y;p)$.  Thus
the ranges of two equivalent paths are equivalent vertices.  As the
next proposition shows, Condition (P3) of Definition
\ref{def:pathequivalence} is enough to ensure that the sources of
equivalent paths are equivalent.

\begin{proposition}\label{welldefprop}
The maps $\widetilde{r}$ and $\widetilde{s}$ are well defined.
\end{proposition}

\begin{proof}
Suppose $(x;(m,n))\sim(y;(p,q))$.   Then (P1) of Definition
\ref{def:pathequivalence} implies that $n\wedge d(x)-m\wedge
d(x)=q\wedge d(y)-p\wedge d(y)$.  Subtracting this from the equation
in (P3) gives
\begin{gather*}
n-m+m\wedge d(x)-n\wedge d(x) = q-p+p\wedge d(y) -q\wedge d(y)\\
\Leftrightarrow n-n\wedge d(x) -(m-m\wedge d(x))=q-q\wedge d(y)-(p-p\wedge d(y))\\
\Leftrightarrow n-n\wedge d(x)=q-q\wedge d(y) \quad \text{using
(P2)}.
\end{gather*}
Since (P1) gives $x(n\wedge d(x))=y(q\wedge d(y))$, it follows  that
$(x;n)\approx(y;q)$.  Therefore, $\widetilde{s}$ is well defined.

To show $\widetilde{r}$ is well-defined, first consider  the case
where $m\leq d(x)$.  Then $m\wedge d(x)=m$.  Therefore, $m-m\wedge
d(x)=0$, and (P2) implies that $p\wedge d(y)=p$.  Hence, $x(m)=y(p)$
by (P1).

If $m\not\leq d(x)$, then (P1) of Definition
\ref{def:pathequivalence}  implies that $x(m\wedge d(x))=y(p\wedge
d(y))$ and thus condition (V1) of Definition
\ref{def:vertexequivalence} is satisfied.  Condition (P2) of
Definition \ref{def:pathequivalence} is precisely (V2) of Definition
\ref{def:vertexequivalence}.  Therefore, $(x;m)\approx (y;p)$, and
$\widetilde{r}$ is well defined.
\end{proof}

\begin{proposition}\label{prop:shiftequivalence}
Suppose $x,y,\in\Linf$, and suppose $p,q\in\N^k$ are such  that $p\leq
d(x)$, $q\leq d(y)$ and $\sigma^p x=\sigma^q y$. For all
$a,b\in\N^k$, if $a\leq b$ and $b+p\not\leq d(x)$, then $b+q\not\leq
d(y)$ and $\pec(x,{a+p},{b+p})=\pec(y,{a+q},{b+q})$.
\end{proposition}

\begin{proof}
By definition, $d(\shift(p,x))=d(x)-p$ and $d(\shift(q,y))=d(y)-q$.
Therefore
\begin{equation}\label{eq:shiftdegrees}
d(x)=d(\shift(p,x))+p \text{ and } d(y)=d(\shift(q,y))+q.
\end{equation}
Suppose $a,b\in\N^k$ are such that $a\leq b$ and $b+p\not\leq d(x)$.
Then
\begin{align*}
b+p\not\leq d(x)&\Leftrightarrow b+p\not\leq d(\shift(p,x))+p\\
&\Leftrightarrow b\not\leq d(\shift(p,x))\\
&\Leftrightarrow b\not\leq d(\shift(q,y))\\
&\Leftrightarrow b\not\leq d(y)-q\\
&\Leftrightarrow b+q\not\leq d(y).
\end{align*}
Thus $\pec(x,{a+p},{b+p})$ and $\pec(y,{a+q},{b+q})$ are elements in
$\widetilde{P_\Lambda}$.  To show that
$\pec(x,a+p,b+p)=\pec(y,a+q,b+q)$, consider
\begin{align*}
x((a+p)\wedge d(x), (b+p)\wedge d(x))&=x((a+p)\wedge (d(\shift(p,x))+p),(b+p)\wedge (d(\shift(p,x))+p))\\
&=x(a\wedge d(\shift(p,x))+p, b\wedge d(\shift(p,x))+p)\\
&=\shift(p,x)(a\wedge d(\shift(p,x)),b\wedge d(\shift(p,x)))\\
&=\shift(q,y)(a\wedge d(\shift(q,y)),b\wedge d(\shift(q,y)))\\
&=y(a\wedge d(\shift(q,y))+q,b\wedge d(\shift(q,y))+q)\\
&=y((a+q)\wedge (d(\shift(q,y))+q),(b+q)\wedge (d(\shift(q,y))+q))\\
&=y((a+q)\wedge d(y),(b+q)\wedge d(y)).
\end{align*}
Thus condition (P1) of Definition \ref{def:pathequivalence} is
satisfied.  To show condition (P2), we have
\begin{align*}
a+p-(a+p)\wedge d(x)&=a+p-(a+p)\wedge (d(\shift(p,x))+p)\\
&=a+p-(a\wedge d(\shift(p,x))+p)\\
&=a-a\wedge d(\shift(p,x))\\
&=a-a\wedge d(\shift(q,y))\\
&=a+q-(a\wedge d(\shift(q,y))+q)\\
&=a+q-(a+q)\wedge(d(\shift(q,y))+q)\\
&=a+q-(a+q)\wedge d(y).
\end{align*}
Condition (P3) is clear.  Hence, $\pec(x,a+p,b+p)=\pec(y,a+q,b+q)$.
\end{proof}

If $p=0$, then $x=\shift(q,y)$ and Proposition
\ref{prop:shiftequivalence} implies that for all $b\not\leq d(x)$, we have
$\pec(x,a,b)=\pec(y,a+q,b+q)$ .

The following proposition will be used to compose two paths in
$\widetilde{P_\Lambda}$.

\begin{proposition}\label{prop:compdefinition}
Let $\pec(x,m,n),\pec(y,p,q)\in\widetilde{P_\Lambda}$ be  such that
$\vcc(x,n)=\vcc(y,p)$.  Define $z=x(0,n\wedge d(x))\sigma^{p\wedge d(y)}y$.  Then
\begin{enumerate}
\item[(i)] $z\in\Lambda^{\leq\infty}$,
\item[(ii)] $m\wedge d(x)=m\wedge d(z)$ and $n\wedge d(x)=n\wedge d(z)$,
\item[(iii)]$x(m\wedge d(x),n\wedge d(x))=z(m\wedge d(z),n\wedge d(z))$ and $y(p\wedge d(y),q\wedge d(y))=z(n\wedge d(z), (n+q-p)\wedge d(z))$
\end{enumerate}
\end{proposition}

\begin{proof}
\emph{Proof of (i)}: Since $y\in\Lambda^{\leq\infty}$, $z$ belongs
to $\Linf$ by \cite[Lemmas 2.10 and 2.11]{RSY2}.
\\

\emph{Proof of (ii)}: We will show the  equalities $m\wedge
d(x)=m\wedge d(z)$ and $n\wedge d(x)=n\wedge d(z)$ hold on a
coordinate by coordinate basis.  Let $i\in\{1,2,\ldots,k\}$.  Since
$\vcc(x,n)=\vcc(y,p)$, it follows that $n-(n\wedge d(x))=p-(p\wedge
d(y))$.  Therefore,
\begin{align*}
d(z)&=(n\wedge d(x))+d(y)-(p\wedge d(y))\\
&=d(y)+n-p.
\end{align*}

Furthermore,  since $n-(n\wedge d(x))=p-(p\wedge d(y))$, $n_i\leq
d(x)_i$ if and only if $p_i\leq d(y)_i$.
\\

\textsc{Case 1:}  Suppose $d(y)_i=\infty$.  Then $p_i< d(y)_i$, and
so $m_i\leq n_i\leq d(x)_i$.  Moreover, $d(z)_i=\infty$ by
definition, so $(m\wedge d(x))_i=m_i=(m\wedge d(z))_i$ and $(n\wedge
d(x))_i=n_i=(n\wedge d(z))_i$.
\\

\textsc{Case 2:} Suppose $d(y)_i<\infty$. We have
$$d(z)_i=d(y)_i +n_i-p_i=d(y)_i+(n\wedge d(x))_i-(p\wedge d(y))_i<\infty.$$
Suppose $p_i\leq d(y)_i$.  Then, as before, $m_i\leq n_i\leq
d(x)_i$.   Also $d(y)_i-p_i\geq 0$.  This implies $m_i\leq n_i\leq
n_i+d(y)_i-p_i=d(z)_i$.  Thus $(m\wedge d(x))_i=m_i=(m\wedge
d(z))_i$ and $(n\wedge d(x))_i=n_i=(n\wedge d(z))_i$.

Next suppose $p_i>d(y)_i$.  Then $n_i>d(x)_i$ as well.  In this case
\begin{align*}
d(z)_i &= (n\wedge d(x))_i+d(y)_i-(p\wedge d(y))_i\\
&= d(x)_i+d(y)_i-d(y)_i\\
&= d(x)_i.
\end{align*}
Consequently $(m\wedge d(x))_i=(m\wedge d(z))_i$ and $(n\wedge
d(x))_i=(n\wedge d(z))_i$.

So in either case, we have $(m\wedge d(x))_i=(m\wedge d(z))_i$  and
$(n\wedge d(x))_i=(n\wedge d(z))_i$.  Since $i$ was arbitrarily
chosen, this proves (ii).
\\

\emph{Proof of (iii)}: Notice that (ii) implies
$$z(m\wedge d(z),n\wedge d(z))=z(m\wedge d(x),n\wedge d(x))=x(m\wedge d(x),n\wedge d(x))$$
because $z=x(0,n\wedge d(x))\sigma^{p\wedge d(y)}y$.   Also
$m-m\wedge d(x)=m- m\wedge d(z)$.  Thus $\pec(x,m,n)=\pec(z,m,n)$.

To show $\pec(z,n,n+q-p)= \pec(y,p,q)$, we have   that
$\shift({n\wedge d(x)},z)=\shift({p\wedge d(y)},y)$.  By (ii), we
have $n\wedge d(z)=n\wedge d(x)$, and since $\vcc(x,n)=\vcc(y,p)$,
it follows that $n-n\wedge d(z)=p-p\wedge d(y)$.  Then
\begin{align*}
&\pec(z,n,{n+q-p})\\
&=\pec(z,{n-n\wedge d(z)+n\wedge d(z)},{n+q-p-n\wedge d(z)+n\wedge d(z)})\\
&=\pec(y,{n-n\wedge d(z)+p\wedge d(y)},{n+q-p-n\wedge d(z)+p\wedge d(y)})\text{ by Proposition \ref{prop:shiftequivalence}}\\
&=\pec(y,{p-p\wedge d(y)+p\wedge d(y)},{p+q-p-p\wedge d(y)+p\wedge d(y)})\\
&=\pec(y,p,q).
\end{align*}
This proves (iii).

\end{proof}

\begin{remark} Condition~(P1) of Definition~\ref{def:pathequivalence} and Proposition~\ref{prop:compdefinition} imply that
$$x(m\wedge d(x),n\wedge d(x))=z(m\wedge d(z),n\wedge d(z)),\text {and}$$
$$y(p\wedge d(y),q\wedge d(y))=z(n\wedge d(z),(n+q-p)\wedge d(z)).$$
\end{remark}

\begin{definition}\label{def:composition}
Let
$\widetilde{P_\Lambda}\times_{\widetilde{V_\Lambda}}\widetilde{P_\Lambda}$
be the set
$$\{(\pec(x,m,n),\pec(y,p,q))\in\widetilde{P_\Lambda}\times\widetilde{P_\Lambda}:
\widetilde{s}(\pec(x,m,n))=\widetilde{r}(\pec(y,p,q))\}.$$ For
$(\pec(x,m,n),\pec(y,p,q)\in\widetilde{P_\Lambda}\times_{\widetilde{V_\Lambda}}\widetilde{P_\Lambda}$,
let $z=x(0,n\wedge d(x))\sigma^{p\wedge d(y)}y$. Define
$$\pec(x,m,n)\circ\pec(y,p,q)=\pec(z,m,n+q-p).$$
\end{definition}

\begin{proposition}\label{prop:compositionwelldefined}
The composition defined on
$\widetilde{P_\Lambda}\times_{\widetilde{V_\Lambda}}\widetilde{P_\Lambda}$
given in Definition \ref{def:composition} is well-defined.
\end{proposition}

\begin{proof}
This follows from Proposition \ref{prop:compdefinition}.
\end{proof}

\begin{remark}\label{compproprmk}
If $\pec(x,m,n), \pec(y,p,q)$, and $\pec(z,m,n+q-p)$ are as above,
notice that $\pec(z,m,n)\circ\pec(z,n,n+q-p)=\pec(z,m,n+q-p)$ as
well.  Thus Proposition \ref{prop:compositionwelldefined} implies
that
$\widetilde{r}(\pec(z,m,n+q-p))=\widetilde{r}(\pec(z,m,n))=\widetilde{r}(\pec(x,m,n))$
and
$\widetilde{s}(\pec(z,m,n+q-p))=\widetilde{s}(\pec(z,n,n+q-p))=\widetilde{s}(\pec(y,p,q))$
by Proposition \ref{welldefprop}.
\end{remark}

\begin{proposition}\label{prop:bicatcompdefinition}
For $\lambda\in\Lambda$ and $(x;(m,n))\in P_\Lambda$ with
$s(\lambda)=x(m)$,  let $z=\lambda\sigma^{m}x$.  Then
\begin{enumerate}
\item[(i)] $z\in\Lambda^{\le\infty}$, \text{ and}
\item[(ii)] $\pec(z,{d(\lambda)},{n-m+d(\lambda)})=\pec(x,m,n)$.
\end{enumerate}
\end{proposition}

\begin{proof}
Since $x\in\Lambda^{\leq\infty}$, (i) follows from \cite[Lemmas 2.10 and 2.11]{RSY2}.

Using the fact that $\shift({d(\lambda)},z)=\shift(m,x)$,
Proposition \ref{prop:shiftequivalence}  implies that
$$\pec(z,{d(\lambda)},{n-m+d(\lambda)})=\pec(x,m,n-m+m)=\pec(x,m,n).$$  Thus (ii) follows.
\end{proof}

\begin{definition}\label{def:bicatcomposition}
Let $\Lambda\times_{\Lambda^0}\widetilde{P_\Lambda}
=\{(\lambda,\pec(x,m,n))\in\Lambda\times\widetilde{P_\Lambda}:s(\lambda)=\widetilde{r}(\pec(x,m,n))\}$.
For
$(\lambda,\pec(x,m,n))\in\Lambda\times_{\Lambda^0}\widetilde{P_\Lambda}$,
let $z=\lambda\shift(m,x)$.   Define
$$\lambda\circ\pec(x,m,n)=\pec(z,0,{d(\lambda)+n-m}).$$
\end{definition}

The proof of the following is a direct consequence of Proposition \ref{prop:bicatcompdefinition}.
\begin{proposition}\label{prop:bicatcompwelldefined}
The composition defined  on
$\Lambda\times_{\Lambda^0}\widetilde{P_\Lambda}$ given in Definition
\ref{def:bicatcomposition} is well-defined.
\end{proposition}

\begin{remark}\label{bicatcompproprmk}
As in Remark \ref{compproprmk}, if $\lambda, \pec(x,m,n)$,  and
$\pec(z,d(\lambda),{n-m+d(\lambda)})$ are as above, Proposition
\ref{prop:bicatcompwelldefined} implies that
$\widetilde{r}(\pec(z,d(\lambda),{n-m+d(\lambda)}))=r(\lambda)$ and
that
$\widetilde{s}(\pec(z,d(\lambda),{n-m+d(\lambda)}))=\widetilde{s}(\pec(x,m,n))$.
\end{remark}

We are now ready to define the $k$-graph $\overline{\Lambda}$
mentioned in Theorem \ref{thm:nosources}.  The objects of
$\overline{\Lambda}$ consist of the objects of $\Lambda$ together
with the elements of $\widetilde{V_\Lambda}$.  The morphisms of
$\overline{\Lambda}$ are the morphisms of $\Lambda$ and the elements
of $\widetilde{P_\Lambda}$.  Definitions \ref{def:composition} and
\ref{def:bicatcomposition} describe the composition
in $\overline{\Lambda}$.

\begin{definition}\label{def:Lambdabar}
Define $\overline{\Lambda}$ by
\begin{gather*}
\Obj(\overline\Lambda)=\Obj(\Lambda)\cup\widetilde{V_\Lambda}\\
\Hom(\overline\Lambda)=\Hom(\Lambda)\cup\widetilde{P_\Lambda},
\end{gather*}
with $\overline{r}$ and $\overline{s}$ defined as follows:
\begin{gather*}
\overline{r}:\Hom(\overline{\Lambda})\to\Obj(\overline{\Lambda})\\
\overline{r}\mid_{\Hom(\Lambda)}=r,\ \text{and
$\overline{r}\mid_{\widetilde{P_\Lambda}}=\widetilde{r}$}
\end{gather*}
and
\begin{gather*}
\overline{s}:\Hom(\overline{\Lambda})\to\Obj(\overline{\Lambda})\\
\overline{s}\mid_{\Hom(\Lambda)}=s, \ \text{and
$\overline{s}\mid_{\widetilde{P_\Lambda}}=\widetilde{s}$}
\end{gather*}

Let
$$\Hom(\overline{\Lambda})\times_{\Obj(\overline{\Lambda})}\Hom(\overline{\Lambda})=(\Lambda\times_{\Lambda^0}\Lambda)
\bigcup(\Lambda\times_{\Lambda^0}\widetilde{P_\Lambda})\bigcup(\widetilde{P_\Lambda}\times_{\widetilde{V_\Lambda}}\widetilde{P_\Lambda}).$$
Define
$\circ:\Hom(\overline{\Lambda})\times_{\Obj(\overline{\Lambda})}\Hom(\overline{\Lambda})\to\Hom(\overline{\Lambda})$
as follows. For
$(\lambda,\pec(x,m,n))\in\Lambda\times_{\Lambda^0}\widetilde{P_\Lambda}$
define
$$\lambda \circ \pec(x,m,n)=\pec(\lambda\sigma^{m}x,0,{d(\lambda)+n-m}).$$

For
$(\pec(x,m,n),\pec(y,p,q))\in\widetilde{P_\Lambda}\times_{\widetilde{V_\Lambda}}\widetilde{P_\Lambda}$,
let $z=x(0,n\wedge d(x))\shift({p\wedge d(y)},y)$ and define
$$\pec(x,m,n)\circ\pec(y,p,q)=\pec(z,m,n+q-p)$$

For $\lambda,\mu\in\Lambda$ define $\lambda\circ\mu$ as in
$\Lambda$.

Define $\id_{\vcc(x,m)}=\pec(x,m,m)$ for
$\vcc(x,m)\in\widetilde{V_\Lambda}$, and  define $\id_v$ as in
$\Lambda$ for $v\in\Obj(\Lambda)$.
\end{definition}

\begin{lemma}\label{lambdabariscategory}
With the definitions given above, $\overline{\Lambda}$ is a
category.
\end{lemma}

\begin{proof}
Using the axioms for a category detailed in \cite[Section I.2]{MacLane}, it must be shown that:

\begin{itemize}
\item[(i)]\label{catcond1} $\overline{r}(\id_{c})=c=\overline{s}(\id_c)$ for all $c\in\Obj(\overline{\Lambda})$;
\item[(ii)]\label{catcond2}$\overline{s}(f\circ g)=\overline{s}(g)$ and $\overline{r}(f\circ g)=\overline{r}(f)$ for all $f,g\in\Hom(\overline{\Lambda})$;
\item[(iii)] \label{catcond3}$(f\circ g)\circ h=f\circ (g\circ h)$ for all $f,g,h\in\Hom(\overline{\Lambda})$;
\item[(iv)] \label{catcond4}$f\circ\id_c=f$ and $\id_c\circ g=g$ for all $c\in\Obj(\overline{\Lambda})$ and $f,g\in\Hom(\overline{\Lambda})$ such that $\overline{s}(f)=c=\overline{r}(g)$.
\end{itemize}

\emph{Proof of (i)}: Since $\overline{r}=r$ and $\overline{s}=s$ on
$\Lambda$, (i) holds  for $v\in\Obj(\Lambda)$ because $\Lambda$ is a
category.  If $\vcc(x,m)\in\widetilde{V_\Lambda}$, then $m\not\leq
d(x)$.  Therefore,
\begin{align*}
\overline{r}(\id_{\vcc(x,m)})&=\overline{r}(\pec(x,m,m))\\
&=\vcc(x,m)\\
&=\overline{s}(\pec(x,m,m))\\
&=\overline{s}(\id_{\vcc(x,m)}).
\end{align*}
Thus (i) is true for all $c\in\Obj(\overline{\Lambda})$.
\\

\emph{Proof of (ii)}: Suppose $\lambda,\mu\in
\Hom(\Lambda)\subseteq\Hom(\overline{\Lambda})$.   Then (ii) follows
because $\Lambda$ is a category and $\overline{s}$ agrees with $s$
on $\Hom(\Lambda)$.  If $\lambda\in\Hom(\Lambda)$ and
$\pec(x,m,n)\in\widetilde{P_\Lambda}\subseteq\Hom(\overline{\Lambda})$,
then $\overline{s}(\pec({\lambda\shift({m\wedge
d(x)},x)},0,{d(\lambda)+n-m}))=\vcc({\lambda\shift({m\wedge
d(x)},x)},{d(\lambda)+n-m})$.  Thus (ii) is true because
$\vcc({\lambda\shift(m,x)},{d(\lambda)+n-m})=\vcc({\shift(m,x)},n-m)=\vcc(x,n)$
by Proposition \ref{prop:shiftequivalence} (applied twice).  To show
that (ii) holds for
$\pec(x,m,n),\pec(y,p,q)\in\widetilde{P_\Lambda}$, the definition of composition in $\overline{\Lambda}$ yields
$\pec(x,m,n)\circ\pec(y,p,q)=\pec({x(0,n\wedge d(x))\shift({p\wedge
d(y)},y)},m,n+q-p)$.  Therefore,
$$\overline{s}(\pec({x(0,n\wedge d(x))\shift({p\wedge d(y)},y)},m,n+q-p))=\vcc({x(0,n\wedge d(x))\shift({p\wedge d(y)},y)},n+q-p),$$
and
\begin{align*}
\vcc({x(0,n\wedge d(x))\shift({p\wedge d(y)},y)},n+q-p)&=\vcc({\shift({p\wedge d(y)},y)},{n+q-p-n\wedge d(x)})\\
&=\vcc(y,{n+q-p-n\wedge d(x)+p\wedge d(y)})\\
&\hspace{.75in}\text{ by Proposition \ref{prop:shiftequivalence}}\\
&=\vcc(y,{n-n\wedge d(x)+q-(p-p\wedge d(x))})\\
&=\vcc(y,q)
\end{align*}
since $\vcc(x,n)=\vcc(y,p)$ implies that $n-n\wedge d(x)=p-p\wedge
d(y)$.

Showing that $\overline{r}(f\circ g)=\overline{r}(f)$ follows in a
similar manner since.
\\

\emph{Proof of (iii)}: There are four cases to consider.
\\

\textsc{Case 1:} Suppose
$\lambda,\mu,\nu\in\Hom(\Lambda)\subseteq\Hom(\overline{\Lambda})$.
Condition (iii) holds in this case because $\Lambda$ is a category
and composition in $\overline{\Lambda}$ on
$\Hom(\Lambda)\subseteq\Hom(\overline{\Lambda})$ agrees with the
composition in $\Lambda$.
\\

\textsc{Case 2:}  Suppose $\lambda, \mu\in\Hom(\Lambda)$ and
$\pec(x,m,n)\in\widetilde{P_\Lambda}\subseteq\Hom(\overline{\Lambda})$.
Then
\begin{align*}
(\lambda\circ\mu)\circ\pec(x,m,n)&=(\lambda\mu)\circ\pec(x,m,n)\\
&=\pec((\lambda\mu)\sigma^mx,0,n-m+d(\lambda\mu))\ \\
&=\pec(\lambda(\mu\sigma^mx),0,{n-m+d(\lambda)+d(\mu)})\\
&\hspace{.75in}\text{because composition in $\Lambda$ is associative}\\
&=\lambda\circ\pec({\mu\sigma^mx},0,{n-m+d(\mu)})\\
&=\lambda\circ(\mu\circ\pec(x,m,n)).
\end{align*}
\\

\textsc{Case 3:} Suppose $\lambda\in\Hom(\Lambda)$ and
$\pec(x,m,n),\pec(y,p,q)\in\widetilde{P_\Lambda}$.  Then
\begin{align*}
(\lambda\circ\pec(x,m,n))\circ\pec(y,p,q)&=\pec(\lambda\sigma^mx,0,{n-m+d(\lambda)})\circ\pec(y,p,q)\\
&=\pec(z,0,{n-m+d(\lambda)+q-p})
\end{align*}
where $z=(\lambda\sigma^mx)(0,(n-m+d(\lambda))\wedge
d(\lambda\sigma^mx))\sigma^{p\wedge d(y)}y$.

On the other hand,
\begin{align*}
\lambda\circ(\pec(x,m,n)\circ\pec(y,p,q))&=\lambda\circ\pec(w,m,n+q-p)\\
&=\pec(\lambda\sigma^mw,0,{n-m+q-p+d(\lambda)})
\end{align*}
where $w=x(0,n\wedge d(x))\sigma^{p\wedge d(y)}y$.

To show that
$$\pec(z,0,n{-m+d(\lambda)+q-p})=\pec(\lambda\sigma^mw,0,{n-m+q-p+d(\lambda)}),$$
notice that

\begin{align}
&z\left(0\wedge d(z), (n-m+q-p+d(\lambda))\wedge d(z)\right)\label{eq:compositioncase3_1}\\
&=z\left(0,(n-m+d(\lambda))\wedge d(z)\right))\notag\circ\\
&\quad\quad\circ z\left((n-m+d(\lambda))\wedge d(z),(n-m+d(\lambda)+q-p)\wedge d(z)\right)\notag\\
&=(\lambda\sigma^mx)\left(0,(n-m+d(\lambda)\right)\wedge
d(\lambda\sigma^mx))y\left(p\wedge d(y),q\wedge d(y)\right)\notag
\end{align}

by Proposition~\ref{prop:compdefinition}~(iii).

Since $d(\lambda\sigma^{m}x)=d(\lambda)-m+d(x)$, we have that
\begin{align*}
n-m+d(\lambda)\wedge d(\lambda\shift(m,x))&=(n-m+d(\lambda)\wedge(d(x)-m+d(\lambda))\\
&=d(\lambda)-m+n\wedge d(x)
\end{align*}
because addition in $\N^k$ distributes over $\wedge$.  Thus we can
continue with the calculation:
\begin{align}
(\lambda\sigma^m&x)(0,(n-m+d(\lambda)\wedge d(\lambda\sigma^mx))y(p\wedge d(y),q\wedge d(y))\label{eq:compositioncase3_2}\\
&=(\lambda\sigma^mx)(0,d(\lambda)-m +n\wedge d(x))y(p\wedge d(y),q\wedge d(y))\notag\\
&=\lambda(\sigma^mx)(0,-m+n\wedge d(x))y(p\wedge d(y),q\wedge d(y))\notag\\
&=\lambda x(m,n\wedge d(x))y(p\wedge d(y),q\wedge d(y)).\notag
\end{align}
Equations~\eqref{eq:compositioncase3_1} and \eqref{eq:compositioncase3_2} show that
\begin{equation}\label{eq:compositioncase3_3}
z\left(0\wedge d(z), (n-m+q-p+d(\lambda))\wedge d(z)\right)=\lambda x(m,n\wedge d(x))y(p\wedge d(y),q\wedge d(y)).
\end{equation}
Similarly, it can be shown using
Proposition~\ref{prop:compdefinition}  that

\begin{equation}\label{eq:compositioncase3_4}
(\lambda\sigma^mw)(0,(n-m+q-p+d(\lambda))\wedge d(\lambda\sigma^mw))=\lambda x(m,n\wedge d(x))y(p\wedge d(y),q\wedge d(y)).\end{equation}

Equations~\eqref{eq:compositioncase3_3} and \eqref{eq:compositioncase3_4} show that Conditions~(P1) of Definition \ref{def:pathequivalence} is satisfied.  Condition~(P2)  holds by Proposition~\ref{prop:compdefinition}~(ii).  Clearly, Condition (P3) holds; therefore
$\pec(z,0,{n-m+d(\lambda)+q-p})=\pec(\lambda\sigma^mw,0,{n-m+q-p+d(\lambda)})$
in $\widetilde{P_\Lambda}$.
\\

\textsc{Case 4:}  Suppose
$\pec(x,m,n),\pec(y,p,q),\pec(z,t,u)\in\widetilde{P_\Lambda}\subseteq\Hom(\overline{\Lambda})$.
We must show that
$$(\pec(x,m,n)\circ\pec(y,p,q))\circ\pec(z,t,u)=\pec(x,m,n)\circ(\pec(y,p,q)\circ\pec(z,t,u)).$$
Let $W_1=x(0,n\wedge d(x)\shift({p\wedge d(y)},y)$.   Then
$\pec(x,m,n)\circ\pec(y,p,q)=\pec(W_1,m,{q-p+n})$.  Next,
define $Z_1=W_1(0,(q-p+n)\wedge d(W_1))\shift({t\wedge d(z)},z)$.
Then,
\begin{align*}
(\pec(x,m,n)\circ\pec(y,p,q))\circ(z,t,u)&=\pec(W_1,m,q-p+n)\circ\pec(z,t,u)\\
&=\pec(Z_1,m,u-t+q-p+n).
\end{align*}
Using Proposition~\ref{prop:compdefinition} (ii) again, we see that
\begin{gather*}
m\wedge d(Z_1)=m\wedge d(W_1)=m\wedge d(x);\\
(q-p+n)\wedge d(Z_1)=(q-p+n)\wedge d(W_1),\text{ and}\\
\pec(Z_1,q-p+n,u-t+q-p+n)=\pec(z,t,u).
\end{gather*}
We then compute
{\allowdisplaybreaks
\begin{align}
&Z_1(m\wedge d(Z_1),(u-t+q-p+n)\wedge d(Z_1))\label{eq:compositioncase4_1}\\
&=Z_1(m\wedge d(Z_1),(q+p-n)\wedge d(Z_1))\circ\notag\\
&\quad \quad\circ Z_1((q+p-n)\wedge d(Z_1),(u-t+q-p+n)\wedge d(Z_1))\notag\\
&=W_1(m\wedge d(W_1),(q-p+n)\wedge d(W_1))z(t\wedge d(z),u\wedge d(z))\text{\quad by Proposition~\ref{prop:compdefinition}~(iii)}\notag\\
&=W_1(m\wedge d(W_1),n\wedge d(W_1))W_1(n\wedge d(W_1),(q+p-n)\wedge d(W_1))z(t\wedge d(z),u\wedge d(z))\notag\\
&=x(m\wedge d(x),n\wedge d(x))y(p\wedge d(y),q\wedge d(y))z(t\wedge
d(z),u\wedge d(z))\text{\quad by Proposition~\ref{prop:compdefinition}~(iii)}.\notag
\end{align}}
Now let $W_2=y(0,q\wedge d(y))\shift({t\wedge d(z)},z)$,  and
$Z_2=x(0,n\wedge d(x))\shift({p\wedge d(W_2)},W_2)$.  From the
definition of composition, it follows that
\begin{align*}
\pec(x,m,n)\circ(\pec(y,p,q)\circ\pec(z,t,u))&=\pec(x,m,n)\circ\pec(W_2,p,u-t+q)\\
&=\pec(Z_2,m,u-t+q-p+n).
\end{align*}
We must show $\pec(Z_1,m,u-t+q-p+n)=\pec(Z_2,m,u-t+q-p+n)$.  An argument similar to that used in Equation~\eqref{eq:compositioncase4_1} proves that
\begin{align}
Z_1(m\wedge d(Z_1),&(u-t+q-p+n)\wedge d(Z_1))\label{eq:compositioncase4_2}\\
&=x(m\wedge d(x),n\wedge d(x))y(p\wedge d(y),q\wedge d(y))z(t\wedge
d(z),u\wedge d(z)).\notag
\end{align}

Equations~\eqref{eq:compositioncase4_1} and \eqref{eq:compositioncase4_2} show that Condition (P1) of
Definition \ref{def:pathequivalence} are satisfied.  Again, Proposition~\ref{prop:compdefinition} (ii) shows Condition (P2) holds.  It is clear
that Condition (P3) is satisfied.  Therefore the equation
$$\pec(Z_1,m,u-t+q-p+n)=\pec(Z_2,m,u-t+q-p+n),$$
holds, and composition is associative in this case.

Hence, $\overline{\Lambda}$ satisfies (iii).
\\

\emph{Proof of (iv)}: Suppose $v$ is an element of $\Lambda^0$
(which is identified with
$\Obj(\Lambda)\subseteq\Obj(\overline{\Lambda})$).  Then (iv)
follows for all
$f,g\in\Hom(\Lambda)$ such that
$s(f)=v=r(g)$ because $\Lambda$ is a category.  There does not exist
any $f\in\widetilde{P_\Lambda}$ such that $\overline{s}(f)=v$.
Suppose $g\in\widetilde{P_\Lambda}$ is such that
$\overline{r}(g)=v$. Then $g=\pec(x,m,n)$ for some $x\in \Linf$ with $x(m)=v$.  Therefore,
\begin{align*}
v\circ \pec(x,m,n)&=\pec({v\shift(m,x)},0,{n-m})\\
&=\pec({\shift(m,x)},0,{n-m})\\
&=\pec(x,m,n) \text{ by Proposition \ref{prop:shiftequivalence}.}
\end{align*}

Next suppose
$\vcc(x,m)\in\widetilde{V_\Lambda}\subseteq\Obj(\overline{\Lambda})$.
There does not exist any
$f\in\Hom(\Lambda)$ such that
$\overline{r}(f)=\vcc(x,m)$ or $\overline{s}(f)=\vcc(x,m)$. Thus, if
$f\in\Hom(\overline{\Lambda})$ such that
$\overline{s}(f)=\vcc(x,m)$, then $f=\pec(y,p,q)$ for some
$\pec(y,p,q)\in\widetilde{P_\Lambda}$ such that
$\vcc(x,m)=\vcc(y,q)$.  Then, by definition of composition in
$\overline{\Lambda}$, we have
\begin{align*}
\pec(y,p,q)\circ\id_{\vcc(x,m)}&=\pec(y,p,q)\circ\pec(x,m,m)\\
&=\pec({y(0,q\wedge d(y))\shift({m\wedge d(x)},x)},p,{q+m-m})\\
&=\pec({y(0,q\wedge d(y))\shift({m\wedge d(x)},x)},p,{q})\\
&=\pec(y,p,q)\quad\text{by Proposition~\ref{prop:compdefinition} (iii).}
\end{align*}
It is shown similarly that  if $\pec(z,t,u)$ is an element of
$\widetilde{P_\Lambda}$ with $\vcc(x,m)=\vcc(z,t)$, the equality
$\id_{\vcc(x,m)}\circ\pec(z,t,u)=\pec(z,t,u)$ holds. We have shown that
(iv) holds, and thus $\overline{\Lambda}$ is a category.
\end{proof}

From now on, we will write $\lambda\mu$ instead of $\lambda\circ\mu$
for  all $\lambda,\mu\in\Hom(\overline{\Lambda})$.

We will view $\N^k$ as a category with one object ($\star$), a morphism set equal to $\N^k$ and with composition determined by addition in $\N^k$.

\begin{definition}\label{def:dbar}
Define $\overline{d}:\overline{\Lambda}\to\N^k$ as follows.   For
all $c\in\Obj(\overline{\Lambda})$, let $\overline{d}(c)=\star$.
Furthermore, define
$$\overline{d}\mid_{\Hom(\Lambda)}=d, \text{ and }\overline{d}(\pec(x,m,n))=n-m, \text{ for } \pec(x,m,n)\in\widetilde{P_\Lambda}.$$
\end{definition}

It is straightforward to show that $\overline{d}$ defines a functor.

\begin{lemma}\label{UFP}
The category $\overline{\Lambda}$ with the functor $\overline{d}$ defined in Definition~\ref{def:dbar}
satisfies  the factorization property.  That is, for
$f\in\Hom(\overline{\Lambda})$ with $\overline{d}(f)=a+b$, there
exist unique elements $g,h\in\Hom(\overline{\Lambda})$ such that
$f=g\circ h$ with $\overline{d}(g)=a$ and $\overline{d}(h)=b$.
\end{lemma}

\begin{proof}
If $f\in\Hom(\Lambda)\subseteq\Hom(\overline{\Lambda})$, then since
$\Lambda$ has the  factorization property and $\overline{d}$ agrees
with $d$ on $\Hom(\Lambda)$, the required elements exist and are
unique.

Suppose that
$\pec(x,m,n)\in\widetilde{P_\Lambda}\subseteq\Hom(\overline{\Lambda})$.
Then $\overline{d}(\pec(x,m,n))=n-m$.  Suppose that $n-m=a+b$. There
are three cases to consider: $m\not\leq d(x)$; $m\leq d(x)$ while
$m+a\not\leq d(x)$; and $m\leq m+a\leq d(x)$.
\\

\textsc{Case 1:} Suppose $m\not\leq d(x)$. By definition of
composition  in $\widetilde{P_\Lambda}$ and Remark
\ref{compproprmk}, the necessary elements exist, namely
$\pec(x,m,m+a)$ and $\pec(x,m+a,n)$.  For uniqueness, suppose that
$\pec(x,m,n)=\pec(x,m,m+a)\pec(x,m+a,n)$ as well as
$\pec(x,m,n)=\pec(y,p,q)\pec(z,t,u)$ with $q-p=a$ and $u-t=b$. Using
the definition of composition in $\widetilde{P_\Lambda}$,
$\pec(y,p,q)\pec(z,t,u)=\pec(w,p,q+u-t)$ where $w=y(0,q\wedge
d(y))\sigma^{t\wedge d(z)}z$.  Since $\overline{\Lambda}$ is a category, it follows that
\begin{gather*}
\vcc(x,m)=\overline{r}(\pec(x,m,n))=\overline{r}(\pec(y,p,q))=\vcc(y,p) \text{ and}\\
\vcc(x,n)=\overline{s}(\pec(x,m,n))=\overline{s}(\pec(z,t,u))=\vcc(z,u).
\end{gather*}
Also, since $\overline{s}(\pec(y,p,q))=\overline{r}(\pec(z,t,u))$,
it  follows that $\vcc(y,q)=\vcc(z,t)$.  Therefore, Condition~(V2) of Definition~\ref{def:vertexequivalence} gives the following equalities:
\begin{equation}\label{eq1UFP}
m-m\wedge d(x)=p-p\wedge d(y)
\end{equation}
\begin{equation}\label{eq2UFP}
n-n\wedge d(x)=u-u\wedge d(z)
\end{equation}
\begin{equation}\label{eq3UFP}
q-q\wedge d(y)=t-t\wedge d(z)
\end{equation}

Furthermore, since
$\pec(x,m,n)=\pec(x,m,m+a)\pec(x,m+a,n)=\pec(y,p,q)\pec(z,t,u)$,  Condition (P1) of Definition~\ref{def:pathequivalence} implies that
\begin{align*}
x(m\wedge d(x), n\wedge d(x))&=x(m\wedge d(x),(m+a)\wedge d(x))x((m+a)\wedge d(x),n\wedge d(x))\\
&=y(p\wedge d(y),q\wedge d(y))z(t\wedge d(z), u\wedge d(z)).
\end{align*}
The first equality above shows that
\begin{equation}\label{eq4UFP}
n\wedge d(x)-m\wedge d(x)=q\wedge d(y)-p\wedge d(y)+u\wedge
d(z)-t\wedge d(z).
\end{equation}

Now, if $q\wedge d(y)-p\wedge d(y)=(m+a)\wedge d(x)-m\wedge
d(x)$, then the  factorization property of $\Lambda$ will imply that
$x(m\wedge d(x),(m+a)\wedge d(x))=y(p\wedge d(y),q\wedge d(y))$.
Then by \eqref{eq1UFP} and the fact that $a=(m+a)-m=q-p$, it will
follow that $\pec(x,m,m+a)=\pec(y,p,q)$.  Consequently, we will have
$\pec(x,m+a,n)=\pec(z,t,u)$.  We will show that $q\wedge
d(y)-p\wedge d(y)=(m+a)\wedge d(x)-m\wedge d(x)$ on a coordinate by
coordinate basis; i.e., by showing that, for all $i\in\{1,2,\dots,k\}$, the equality$(q\wedge d(y))_i-(p\wedge d(y))_i=((m+a)\wedge d(x))_i-(m\wedge d(x))_i$ holds.

Fix $i\in\{1,2,\dots,k\}$. Then \eqref{eq1UFP} implies that $m_i\leq
d(x)_i$ if  and only if $p_i\leq d(y)_i$.  Similarly, by
\eqref{eq2UFP}, $n_i\leq d(x)_i$ if and only if $u_i\leq d(z)_i$,
while \eqref{eq3UFP} ensures $q_i\leq d(y)_i$ if and only if
$t_i\leq d(z)_i$.  Therefore there are 5 cases to consider:

\begin{itemize}
\item[(1-i)] $p_i\leq q_i\leq d(y)_i$ and $m_i\leq m_i+a_i\leq d(x)_i$
\item[(1-ii)] $p_i\leq q_i\leq d(y)_i$ and $m_i\leq d(x)_i< m_i+a_i$
\item[(1-iii)] $p_i\leq d(y)_i< q_i$ and $m_i\leq m_i+a_i\leq d(x)_i$
\item[(1-iv)] $p_i\leq d(y)_i<q_i$ and $m_i\leq d(x)_i< m_i+a_i$
\item[(1-v)] $d(y)_i<p_i\leq q_i$ and $d(x)_i< m_i\leq m_i+a_i$.
\end{itemize}

Cases (1-i) and (1-v) are shown by a simple calculation.

For Case (1-iv), since $n_i\geq m_i+a_i>d(x)_i$, it follows that
$u_i>d(z)_i$.   Furthermore, the fact that $q_i>d(y)_i$ gives the
inequality $t_i>d(z)_i$.  Substituting into \eqref{eq4UFP}, we
obtain
$$d(x)_i-m_i=d(y)_i-p_i+d(z)_i-d(z)_i=d(y_i)-p_i,$$
which shows precisely that $((m+a)\wedge d(x))_i-(m\wedge
d(x))_i=(q\wedge d(y))_i-(p\wedge d(y))_i$.

We will show that the remaining cases cannot, in fact, occur. For
Case (1-ii),  since $m_i\leq d(x)_i<m_i+a_i\leq n_i$, we have
\begin{equation}\label{eq1aUFP}
(n\wedge d(x)-m\wedge d(x))_i=d(x)_i-m_i.
\end{equation}
Also $q_i\leq d(y)_i$ implies that $t_i\leq d(z)_i$.  Since
$d(x)_i<n_i$  guarantees $d(z)_i<u_i$, it follows that
\begin{equation}\label{eq2aUFP}
(q\wedge d(y)-p\wedge d(y))_i+(u\wedge d(z)-t\wedge
d(z))_i=q_i-p_i+d(z)_i-t_i=a_i+d(z)_i-t_i.
\end{equation}
Substituting (\ref{eq1aUFP}) and (\ref{eq2aUFP}) into
(\ref{eq4UFP}),  implies that
\begin{align*}
d(x)_i-m_i&=a_i+d(z)_i-t_i\\
\Longleftrightarrow d(x)_i&=m_i+a_i+d(z)_i-t_i\\
&\geq m_i+a_i.
\end{align*}
But this, with the hypothesis of Case (1-ii),   means
$$d(x)_i<m_i+a_i\leq d(x)_i,$$
which is a contradiction.  Thus Case (1-ii) does not occur.

For Case (1-iii),  since $q_i>d(y)_i$, it is the case that
$d(z)_i<t_i\leq u_i$.   Therefore $n_i>d(x)_i$.  Using
(\ref{eq4UFP}) again we have
\begin{align*}
d(x)_i-m_i&=(n\wedge d(x)-m\wedge d(x))_i\\
&=(q\wedge d(y)-p\wedge d(y))_i+(u\wedge d(z)-t\wedge d(z))_i\\
&=d(y)_i-p_i+d(z)_i-d(z)_i\\
&=d(y)_i-p_i.
\end{align*}
Therefore $d(x)_i-m_i-a_i=d(y)_i-p_i-a_i=d(y)_i-q_i$.  However, the
conditions of  Case (1-iii) imply that $d(x)_i-(m_i+a_i)\geq 0$ and
$d(y)_i-q_i<0$, which is a contradiction.  Consequently, Case (1-iii)
does not occur.
\\

\textsc{Case 2:} Suppose that $m\leq d(x)$ and $m+a\not\leq d(x)$.
By definition of  composition in $\widetilde{P_\Lambda}$ and Remark
\ref{compproprmk}, we may write $\pec(x,m,n)=\pec(x,m,m+a)\pec(x,m+a,n)$.  For uniqueness, suppose
that $\pec(x,m,n)=\pec(y,p,q)\pec(z,t,u)$ with $q-p=a$ and $u-t=b$ as well.  As in Case 1, since
$\pec(x,m,n)=\pec(y,p,q)\pec(z,t,u)$, we have the equalities
\begin{gather*}
x(m)=\overline{r}(\pec(x,m,n))=\overline{r}(\pec(y,p,q))=y(p)\text{ and}\\
\vcc(x,n)=\overline{s}(\pec(x,m,n))=\overline{s}(\pec(z,t,u))=\vcc(z,u)
\end{gather*}

Condition (P1)
of Definition~\ref{def:pathequivalence} implies that
\begin{align*}
x(m,n\wedge d(x))&=x(m,(m+a)\wedge d(x))x((m+a)\wedge d(x),n\wedge d(x))\\
&=y(p,q\wedge d(y))z(t\wedge d(z),u\wedge d(z)),
\end{align*}
and therefore in this case, equation~\eqref{eq4UFP} is replaced with
\begin{equation}\label{eq4aUFP}
n\wedge d(x)-m=q\wedge d(y)-p+u\wedge d(z)-t\wedge d(z).
\end{equation}
Since $m\leq d(x)$,  it follows that $p\leq d(y)$, buy equations
\eqref{eq2UFP} and \eqref{eq3UFP} still hold.

The factorization property of $\Lambda$ will give the uniqueness
provided that
$$(m+a)\wedge d(x)-m=q\wedge d(y)-p.$$
Again, this will be done on a coordinate by coordinate basis.  Fix
$i\in\{1,2,\dots,k\}$.   This time there are four cases to consider:
\begin{itemize}
\item[(2-i)] $p_i\leq q_i\leq d(y)_i$ and $m_i\leq m_i+a_i\leq d(x)_i$
\item[(2-ii)] $p_i\leq q_i\leq d(y)_i$ and $m_i\leq d(x)_i< m_i+a_i$
\item[(2-iii)] $p_i\leq d(y)_i< q_i$ and $m_i\leq m_i+a_i\leq d(x)_i$
\item[(2-iv)] $p_i\leq d(y)_i< q_i$ and $m_i\leq d(x)_i< m_i+a_i$.
\end{itemize}
Cases (2-i) is a simple calculation.  The same argument used to
prove Case (1-iv)  proves Case (2-iv).  We will show the remaining
two cases cannot occur.

For Case (2-ii), $d(x)_i<m_i+a_i\leq n_i$.  Therefore $d(z)_i<u_i$.
Since $q_i\leq d(y)_i$, it follows that $t_i\leq d(z)_i$.  Then
using \eqref{eq4aUFP},
\begin{align*}
d(x)_i-m_i&=(n\wedge d(x)-m)_i\\
&=(q\wedge d(y)-p)_i+(u\wedge d(z)-t\wedge d(z))_i\\
&=q_i-p_i+d(z)_i-t_i\\
&=a_i+d(z)_i-t_i\\
&\geq a_i
\end{align*}
\noindent since $d(z)_i-t_i\geq 0$.  This gives that $d(x)_i\geq
m_i+a_i$, a contradiction.

For Case (2-iii), $d(y)_i<q_i$ implies $d(z)_i<t_i\leq u_i$, and
hence that $d(x)_i<n_i$.  Thus
\begin{align*}
d(x)_i-m_i&=(n\wedge d(x)-m)_i\\
&=(q\wedge d(y)-p)_i+(u\wedge d(z)-t\wedge d(z))_i\\
&=d(y)_i-p_i+d(z)_i-d(z)_i\\
&=d(y)_i-p_i.
\end{align*}
Therefore, using the argument in Case 1, it follows that
$d(x)_i-(m_i+a_i)=d(y)_i-q_i$.   This is a contradiction since
$d(x)_i-(m_i+a_i)\geq 0$ and $d(y)_i-q_i<0$.
\\

\textsc{Case 3:}  Suppose that $m\leq d(x)$ and $m+a\leq d(x)$. Then
using the  definition of composition in $\overline{\Lambda}$, we
have $\pec(x,m,n)=x(m,m+a)\pec(x,m+a,n)$.  To show uniqueness,
suppose that $\pec(x,m,n)=\lambda\pec(y,p,q)$ for some
$\lambda\in\Hom(\Lambda)$ and $\pec(y,p,q)\in\widetilde{P_\Lambda}$,
with $\overline{d}(\lambda)=a$ and $q-p=b=n-(m+a)$.  Then by
Condition (P1) of Definition~\ref{def:pathequivalence},
\begin{align*}
x(m,n\wedge d(x))&=x(m,m+a)x(m+a,n\wedge d(x))\\
&=\lambda y(p,q\wedge d(y).
\end{align*}
The factorization property of $\Lambda$ gives that
$x(m,m+a)=\lambda$.   Consequently, the equality $x(m+a,n\wedge
d(x))=y(p,q\wedge d(y)$ holds.  So, $\pec(x,m+a,n)=\pec(y,p,q)$
which gives uniqueness in this case.
\end{proof}

\begin{theorem}\label{Lambdabarproperties}
Let $(\Lambda,d)$ be a $k$-graph.  Then the extension of this $k$-graph given by the pair
$(\overline{\Lambda},\overline{d})$ of  Definition
\ref{def:Lambdabar} is a $k$-graph with no sources.
\end{theorem}

\begin{proof}
The fact that $\overline{\Lambda}$ is a $k$-graph follows from
Lemmas \ref{lambdabariscategory} and \ref{UFP}.

We will show that $v{\overline{\Lambda}}^{e_i}$ is nonempty for  all
$v\in\overline{\Lambda}^0$ and all $i\in\{1,2,\ldots,k\}$.  If
$v\in{\overline{\Lambda}}^0\backslash\Lambda^0$, then $v=\vcc(x,m)$
for some $x\in\Lambda^{\leq\infty}$ and $m\not\leq d(x)$.  Then
$\pec(x,m,m+e_i)\in v\overline{\Lambda}^{e_i}$ for all
$i\in\{1,2,\ldots,k\}$.  If $v\in\Lambda^0$, choose $x\in
v\Lambda^{\leq\infty}$, which is nonempty by \cite[Lemma 2.11]{RSY2}.  Fix $i\in\{1,2,\ldots,k\}$.  If
$d(x)_i>0$, then $x(0,e_i)\in v\Lambda^{e_i}\subseteq
v\overline{\Lambda}^{e_i}$.  If $d(x)_i=0$, then $\pec(x,0,e_i)\in
v\overline{\Lambda}^{e_i}$.  Hence, for all
$v\in\overline{\Lambda}^0$ and $i\in\{1,2,\ldots,k\}$,
$v\overline{\Lambda}^{e_i}\neq\emptyset$.  Therefore,
$\overline{\Lambda}$ is a $k$-graph without sources.

\end{proof}

Notice that Definition~\ref{def:Lambdabar} provides a way to extend \emph{any} $k$-graph to a larger $k$-graph without sources.  We will show next that if $\Lambda$ is finitely aligned or row-finite, then the extension $\overline{\Lambda}$ will have the same property.

\begin{lemma}\label{lem:lambdamin}
Let $(\Lambda,d)$ be a finitely aligned $k$-graph and let $(\overline{\Lambda},\overline{d})$ be the $k$-graph given in Definition~\ref{def:Lambdabar}.  For $\lambda,\mu \in \Lambda$, we have
$\Lmin(\lambda,\mu)=\Lbarmin(\lambda,\mu)$.
\end{lemma}

\begin{proof}
Of course $\Lmin(\lambda,\mu)\subseteq\Lbarmin(\lambda,\mu)$ because
$\Lambda\subseteq\overline{\Lambda}$.   To show the other
containment, suppose there exists
$(\pec(x,m,n),\pec(y,p,q))\in\Lbarmin(\lambda,\mu)\backslash
\Lmin(\lambda,\mu)$.  Then
\begin{equation}\label{pathequal1}
\pec({\lambda\shift(m,x)},0,{n-m+d(\lambda)})=\lambda
\pec(x,m,n)=\mu \pec(y,p,q)=\pec({\mu\shift(p,y)},0,{q-p+d(\mu)})
\end{equation}
where $\overline{d}(\lambda\pec(x,m,n))=d(\lambda)\vee\d(\mu)$.
Therefore $n-m=d(\lambda)\vee d(\mu)-d(\lambda)$.  But
\eqref{pathequal1} and Condition~(P1) of Definition~\ref{def:pathequivalence} imply
$$\lambda x(m,n\wedge d(x))=\mu y(p,q\wedge d(y)).$$
Since both $\lambda$ and $\mu$ are subpaths of $\lambda x(m,n\wedge
d(x)$,  this implies that $d(\lambda)+n\wedge d(x)-m\geq
d(\lambda)\vee d(\mu)$.  Hence,
$$n\wedge d(x)-m\geq d(\lambda)\vee d(\mu)-d(\lambda)= n-m.$$
It follows that $n\wedge d(x)=n$, and so $n\leq d(x)$,
contradicting  our assumption that the path $\pec(x,m,n)$ is not an element of $\Lambda$.  Thus the set $\Lbarmin(\lambda,\mu)$ is a subset of $\Lmin(\lambda,\mu)$, completing the proof.
\end{proof}

Let $\lambda$ and $\mu$ be two paths in a $k$-graph $\Lambda$.  Recall from  Definition~\ref{def:minimalcommonextension},  that if $(\alpha,\beta)$ is an element of $\Lmin(\lambda,\mu)$, then the path $\lambda\alpha=\mu\beta$ is a minimal common extension of $\lambda$ and $\mu$.  We denote the set of all minimal common extensions of $\lambda$ and $\mu$ by $\MCE(\lambda,\mu)$.  Therefore $\Lambda$ is finitely aligned if and only if $\MCE(\lambda,\mu)$ is finite for all $\lambda,\mu\in\Lambda$.

\begin{theorem}\label{thm:finitelyalignedpreservedinLbar}
Let $(\Lambda,d)$ be a $k$-graph and let $(\overline{\Lambda},\overline{d})$ be the $k$-graph given in
Definition~\ref{def:Lambdabar}.  If $\Lambda$ is finitely aligned, the extension, $\overline{\Lambda}$, is also finitely aligned.  If $\Lambda$ is row-finite, then so is $\overline{\Lambda}$.
\end{theorem}

\begin{proof}
Suppose that $\Lambda$ is finitely aligned.  To show $\overline{\Lambda}$ is finitely aligned, we will show that $|\MCE(\lambda,\mu)|<\infty$ for all $\lambda,\mu\in\overline{\Lambda}$.  Fix two paths $\lambda$ and $\mu$ in $\overline{\Lambda}$.  Let $L=\overline{d}(\lambda)\vee\overline{d}(\mu)$.

First, if $\overline{r}(\lambda)\neq\overline{r}(\mu)$, then $\MCE(\lambda,\mu)=\emptyset$.

Next, suppose $\lambda=\pec(x,m,n)$ and $\mu=\pec(y,p,q)$ are elements of $\overline{\Lambda}\backslash\Lambda$ such that $\overline{r}(\lambda)=\overline{r}(\mu)$.  Any element of $\MCE(\lambda,\mu)$ is of the form $\pec(z,a_z,a_z+L)$ for some $a_z\in\N^k$ and $z\in\Linf$ with $a_z+L\not\leq d(z)$.  Furthermore $\pec(z,a_z,{a_z+d(\lambda)})=\lambda$ and $\pec(z,a_z,{a_z+d(\mu)})=\mu$.  Let $\xi_z=z(a_z\wedge d(z),(a_z+L)\wedge d(z))$.  Then $\xi_z\in\Lambda$, and by Proposition~\ref{prop:compdefinition}, we have that
\begin{align*}
\xi_z&=z(a_z\wedge d(z),(a_z+d(\lambda))\wedge d(z))z((a_z+d(\lambda)\wedge d(z),(a_z+L)\wedge d(z))\\
&=x(m\wedge d(x),n\wedge d(x))z((a_z+d(\lambda)\wedge d(z),(a_z+L)\wedge d(z)).
\end{align*}
Also the degree of $x(m\wedge d(x),n\wedge d(x))$ is
\begin{align*}
n\wedge d(x)-m\wedge d(x)&=(n-m)\wedge d(x)\\
&=(n-m)\wedge d(z) \quad\text{by Proposition~\ref{prop:compdefinition}~(ii)}\\
&=d(\lambda)\wedge d(z).
\end{align*}
On the other hand,
$$\xi_z=y(p\wedge d(y),q\wedge d(y))z((a_z+d(\mu)\wedge d(z),(a_z+L)\wedge d(z)),$$
and $q\wedge d(y)-p\wedge d(y)=d(\mu)\wedge d(z)$.  Since $(d(\lambda)\wedge d(z))\vee(d(\mu)\wedge d(z))=(d(\lambda)\vee d(\mu))\wedge d(z)=d(\xi_z)$, it follows that $\xi_z$ is a minimal common extension of $x(m\wedge d(x),n\wedge d(x))$ and $y(p\wedge d(y),q\wedge d(y))$.

Suppose $\pec(z,a_z,a_z+L)$ and $\pec(w,a_w,a_w+L)$ are two distinct elements of $\MCE(\lambda,\mu)$.  It is clear that Condition (P3) of Definition~\ref{def:pathequivalence} is satisfied, and Condition (P2) holds because $\overline{r}(\lambda)=\overline{r}(\mu)$.  Therefore, Condition (P1) is not satisfied.  This implies that $\xi_z$ and $\xi_w$ are two distinct elements of $\MCE(x(m\wedge d(x),n\wedge d(x)),y(p\wedge d(y),q\wedge d(y)))$.  It follows that $|\MCE(\lambda,\mu)|=|\MCE(x(m\wedge d(x),n\wedge d(x)),y(p\wedge d(y),q\wedge d(y)))|$, which is finite because $\Lambda$ is finitely aligned.

If $\lambda\not\in\Lambda$ and $\mu\in\Lambda$ such that $\overline{r}(\lambda)=\overline{r}(\mu)$, then $\lambda$ may be written as $\pec(x,0,n)$ for some $x\in\Linf$ with $x(0)=r(\mu)$.  In this case every element in $\MCE(\lambda,\mu)$ is if the form $\pec(z,0,L)$ and $z(0,L\wedge d(z))\in\MCE(x(0,n\wedge d(x)),\mu)$.  An argument similar to the previous case shows that $|\MCE(\pec(x,0,n),\mu)|=|\MCE(x(0,n\wedge d(x)),\mu)|$ which again is finite because $\Lambda$ is finitely aligned.

For $\lambda,\mu\in\Lambda$, Proposition~\ref{lem:lambdamin} implies that $\Lbarmin(\lambda,\mu)=\Lmin(\lambda,\mu)$. Thus $\Lbarmin(\lambda,\mu)$ is finite since $\Lambda$ is finitely aligned.

Therefore if $\Lambda$ is finitely aligned, $\MCE(\lambda,\mu)$ is finite for all $\lambda,\mu\in\overline{\Lambda}$, showing that $\overline{\Lambda}$ is finitely aligned.

Now, suppose that $\Lambda$ is row-finite; fix $v\in\overline{\Lambda}^0$
and $i\in\{1,2,\ldots,k\}$.   Since
$\Lambda\subseteq\overline{\Lambda}$ is row finite, the set
$v\Lambda^{e_i}$ is at most finite. Let
$P=v\overline{\Lambda}^{e_i}\backslash\Lambda$.  Any element in $P$
is of the form $\pec(x,m_x,m_x+e_i)$ for some boundary path
$x\in\Lambda^{\leq\infty}$ and $m_x\in\N^k$. For $x,m_x$ and $y,m_y$
such that $\pec(x,m_x,{m_x+e_i}),\pec(y,m_y,{m_y+e_i})\in P$, we
have
$\overline{r}(\pec(x,m_x,{m_x+e_i}))=v=\overline{r}(\pec(y,m_y,{m_y+e_i}))$.
Therefore $m_x-m_x\wedge d(x)=m_y-m_y\wedge d(y)$ and
$$x(m_x\wedge d(x))=y(m_y\wedge d(y))=w$$
for some $w\in\Lambda^0$.  Furthermore, all paths in $P$ have
degree $e_i$.   It follows that two paths in $P$ are distinct if and
only if
$$x(m_x\wedge d(x),(m_x+e_i)\wedge d(x))\neq y(m_y\wedge d(y),(m_y+e_i)\wedge d(y)).$$
Hence, $|P|$ is equal to
$$|\{x(m_x\wedge d(x),(m_x+e_i)\wedge d(x)):\vcc(x,{m_x},{m_x+e_i})\in P\}|.$$
Because $\Lambda$ is row-finite and $P$ is a subset of $\{w\}\cup w\Lambda^{e_i}$,  $P$ is a finite set.
Thus the $k$-graph $\overline{\Lambda}$ is row-finite.
\end{proof}

If $\{t_\lambda:\lambda\in\overline{\Lambda}\}$ is a Cuntz-Krieger
$\overline{\Lambda}$-family, we will show that
$\{t_\lambda:\lambda\in\Lambda\}$ is a Cuntz-Krieger
$\Lambda$-family.  The key elements to show this are Lemma~\ref{lem:lambdamin}, which proves that
$\Lbarmin(\lambda,\mu)$ equals $\Lmin(\lambda,\mu)$ for paths
$\lambda,\mu\in\Lambda$, and the following lemma that shows any finite
exhaustive subset $E$ of $\Lambda$ is exhaustive
in $\overline{\Lambda}$.

\begin{lemma}\label{lem:finiteexhaustive}
Let $(\Lambda,d)$ be a finitely aligned $k$-graph and let $(\overline{\Lambda},\overline{d})$ be the $k$-graph given in Definition~\ref{def:Lambdabar}.  Suppose $v\in \Lambda^0$ and $E\subseteq v\Lambda$ is a finite exhaustive  subset of $\Lambda$.  Then $E$ is also finite exhaustive
subset of $\overline{\Lambda}$.
\end{lemma}

\begin{proof}
Since $E$ is a finite exhaustive subset of $\Lambda$, for every
$\lambda \in\Lambda$  such that $r(\lambda)=v$, there exists $\mu\in
E$ with $\Lmin(\lambda,\mu)\not= \emptyset$.  Therefore, it remains
to show the same holds for paths in $v\overline{\Lambda}\backslash\Lambda$.

Fix $\pec(x,m,n)\in \overline{\Lambda}$ with
$\overline{r}(\pec(x,m,n))=x(m)=v$. We may assume, without loss of
generality, that $m=0$ because $\pec(x,m,n)=\pec(\sigma^m x,0,n-m)$
by Proposition \ref{prop:shiftequivalence}.

Since $x\in\Lambda^{\leq\infty}$, by definition there exists $n_x\in
\N^k$ such  that $n_x\leq d(x)$ and such that if $p\in\N^k$, with $n_x\leq
p\leq d(x)$ and $p_i=d(x)_i$, then $x(p)\Lambda^{e_i}=\emptyset$.
Define
\begin{gather*}
\lambda = x(0,(n\wedge d(x))\vee n_x),\\
\xi =x(0,n\wedge d(x)), \text{ and}\\
\eta=x(n\wedge d(x), (n\wedge d(x))\vee n_x).
\end{gather*}
Notice that if $(n\wedge d(x))_i=d(x)_i$ for some $i$, then
$((n\wedge d(x))\vee n_x)_i=d(x)_i$.   This implies that
$d(\eta)_i=0$ and that
$s(\eta)\Lambda^{e_i}=s(\lambda)\Lambda^{e_i}=\emptyset$.  Hence,
\begin{gather*}\tag{$\star$}
\text{for any path $\alpha\in\Lambda$ such that $r(\alpha)=s(\lambda)=s(\eta)$,}\\
\text{$d(\alpha)_i=0$ if $(n\wedge d(x))_i=d(x)_i$.}
\end{gather*}

There exists $\mu\in E$ and $(\alpha,\beta)\in\Lmin(\lambda,\mu)$
because $E$ is a  finite exhaustive subset of $\Lambda$.  Thus
$$\lambda\alpha=\mu\beta \text{ and } d(\lambda\alpha)=d(\lambda)\vee d(\mu)=((n\wedge d(x))\vee n_x)\vee d(\mu).$$
Since $\lambda=\xi\eta$, and $\Lmin(\lambda,\mu)\not=\emptyset$, it
follows  that $\Lmin(\xi,\mu)\not=\emptyset$.  In particular,  let
\begin{gather*}
\nu=(\lambda\alpha)(n\wedge d(x),(n\wedge d(x))\vee d(\mu))=(\eta\alpha)(0,(n\wedge d(x))\vee d(\mu)-n\wedge d(x)) \\
\text {and } \omega=(\lambda\alpha)(d(\mu),(n\wedge d(x))\vee
d(\mu)).
\end{gather*}
Then $(\nu,\omega)\in\Lmin(\xi,\mu)$.   Moreover, if
$i\in\{1,2,\ldots,k\}$ satisfies $(n\wedge d(x))_i=d(x)_i$,  then
$0=d(\eta)_i=d(\alpha)_i$, giving $d(\nu)_i=0$ and
$d(\xi\nu)_i=d(\xi)_i$.

There exists $y\in\Lambda^{\leq\infty}$ such  that
$y(0,d(\eta\alpha))=\eta\alpha$ by \cite[Lemmas 2.10 and 2.11]{RSY2}.  Also, if $(n\wedge d(x))_i=d(x)_i$, then
$d(y)_i=0$ by ($\star$).
\\

\textsc{Claim 1:} Consider $n-n\wedge d(x)$.   We claim that
$n-n\wedge d(x)\not\leq d(y)$.
\\

\emph{Proof of Claim 1:} To see this, note that  because $n\not\leq
d(x)$ there exists $i\in\{1,2,\ldots,k\}$ such that $n_i> d(x)_i\geq
0$.  Thus $n_i\wedge d(x)_i=d(x)_i$, and $y_i=0$ by the previous
paragraph.  Hence $n_i-(n\wedge d(x))_i=n_i-d(x)_i>0=d(y)_i$, giving
$(n-n\wedge d(x))_i\not\leq d(y)_i$.  This proves Claim 1.
\\

Claim 1 establishes that both the vertex $\vcc(y,{n-n\wedge d(x)})$ and
the path $\pec(y,n-n\wedge d(x), {n-n\wedge d(x)+d(\nu)})$ are
elements in $\overline{\Lambda}$.
\\

\textsc{Claim 2:}  The vertices $\vcc(y,{n-n\wedge d(x)})$ and $\vcc(x,n)$ are equal.
\\

\emph{Proof of Claim 2:} If $(n\wedge d(x))_i=n_i$, then $(n-n\wedge
d(x))_i=0$, and $((n-n\wedge d(x))\wedge d(y))_i=0$.  If $(n\wedge
d(x))_i=d(x)_i$, then $d(y)_i=0$, and $((n-n\wedge d(x))\wedge
d(y))_i=0$.  Therefore $(n-n\wedge d(x))\wedge d(y)=0$. It follows
that
$$y((n-n\wedge d(x))\wedge d(y))=y(0)=r(\eta)=x(n\wedge d(x)).$$
Also $n-n\wedge d(x)-((n-n\wedge d(x))\wedge d(y))=n-n\wedge d(x)$,
which  implies that $\vcc(y,{n-n\wedge d(x)})=\vcc(x,n)$.  This
proves Claim 2.
\\

Claim 2 implies  that $\pec(x,0,n)$ and $\pec(y,{n-n\wedge
d(x)},{n-n\wedge d(x)+d(\nu)})$ are composable in
$\overline{\Lambda}$.  Composing them produces
$$\pec(x,0,n)\pec(y,{n-n\wedge d(x)},{n-n\wedge d(x)+d(\nu)})=\pec({x(0,n\wedge d(x))y},0,{n+d(\nu)}).$$
\\

\textsc{Claim 3:}  We  claim $$\overline{d}(\pec({x(0,n\wedge
d(x))y},0,{n+d(\nu)}))=n+d(\nu)=n\vee d(\mu).$$
\\

\emph{Proof of Claim 3:} Since $d(\nu)=(n\wedge d(x))\vee
d(\mu)-n\wedge d(x)$, we have
\begin{align*}
n+d(\nu)&=n+(n\wedge d(x))\vee d(\mu)-n\wedge d(x)\\
&=n-n\wedge d(x)+(n\wedge d(x))\vee d(\mu)\\
&=(n-n\wedge d(x)+n\wedge d(x))\vee (n-n\wedge d(x)+d(\mu)) \text{ (distributing over $\vee$)}\\
&=n\vee (n-n\wedge d(x)+d(\mu)).
\end{align*}

If $(n\wedge d(x))_i=n_i$, then $(n-n\wedge
d(x)+d(\mu))_i=d(\mu)_i$,  so $(n\vee (n-n\wedge
d(x)+d(\mu)))_i=(n\vee d(\mu))_i$.

On the other hand, if $(n\wedge
d(x))_i=d(x)_i$, then $d(\nu)_i=0$.  But since distributing over
$\wedge$ gives
$$d(\nu)=(n\wedge d(x))\vee d(\mu)-n\wedge d(x)=0\vee (d(\mu)-n\wedge d(x)),$$
it follows that
\begin{equation}\label{finiteexhaustiveeq1}
0\geq (d(\mu)-n\wedge d(x))_i=d(\mu)_i-d(x)_i
\end{equation}
and furthermore,
$$n_i\geq d(\mu)_i+n_i-d(x)_i=d(\mu)_i+n_i-(n\wedge d(x))_i.$$
Therefore $(n\vee (n-n\wedge d(x)+d(\mu)))_i=n_i$.  However, \eqref
{finiteexhaustiveeq1}  implies that $d(x)_i\geq d(\mu)_i$.  Since
$(n\wedge d(x))_i=d(x)_i$, it follows that $n_i\geq d(x)_i\geq
d(\mu)_i$.  Thus, $(n\vee d(\mu))_i=n_i$ as well, establishing Claim
3.
\\

Recall that $y(0,d(\eta\alpha))=\eta\alpha$.  This implies
$$x(0,n\wedge d(x))y=\xi y=(\xi\eta\alpha)\sigma^{d(\eta\alpha)}y
=(\lambda\alpha)\sigma^{d(\eta\alpha)}y=(\mu\beta)\sigma^{d(\eta\alpha)}y$$
because $(\alpha,\beta)\in\Lmin(\lambda,\mu)$.  Hence
\begin{align*}
\pec({x(0,n\wedge d(x))y},0,{n+d(\nu)}) &=\pec({\mu\beta\sigma^{d(\eta\alpha)}y},0,{n+d(\nu)})\\
&=\mu\pec(\beta\sigma^{d(\eta\alpha)}y,d(\mu),{n+d(\nu)-d(\mu)}).
\end{align*}
By Claim 2,
$$\pec({x(0,n\wedge d(x))y},0,{n+d(\nu)})=\pec(x,0,n)\pec(y,n-n\wedge d(x),{n-n\wedge d(x)+d(\nu)}).$$
Claim 3 shows that $\overline{d}(\pec({x(0,n\wedge
d(x))y},0,{n+d(\nu)}))=n\vee d(\mu)$.   Therefore, the path $\pec({x(0,n\wedge d(x))y},0,{n+d(\nu)})$ is a minimal common extension of $\pec(x,0,n)$ and $\mu$. The
pair
$$(\pec(y,n-n\wedge d(x),{n-n\wedge d(x)+d(\nu)}), \pec(\beta\sigma^{d(\eta\alpha)}y,d(\mu),{n+d(\nu)-d(\mu)})$$
is an element of $\Lmin({\pec(x,0,n)},\mu)$, showing that $E$ is a
finite  exhaustive subset of $\overline{\Lambda}$.
\end{proof}

The proof of the next theorem follows easily from Lemmas
\ref{lem:lambdamin}  and \ref{lem:finiteexhaustive}.

\begin{theorem}\label{containsCKlambdafamily}
Let $(\Lambda,d)$ be a finitely aligned $k$-graph and let $(\overline{\Lambda},\overline{d})$ be the $k$-graph given in Definition~\ref{def:Lambdabar}.  If $\{t_{\lambda}:\lambda\in\overline{\Lambda}\}$ is a Cuntz-Krieger
$\overline{\Lambda}$-family,  then the restriction of this set to the elements generated by the subgraph $\Lambda$,
$\{t_{\lambda}:\lambda\in\Lambda\}$, is a Cuntz-Krieger
$\Lambda$-family.
\end{theorem}

\begin{proof}
Conditions (TCK1) and (TCK2) of Definition
\ref{def:Toeplitz-Cuntz-Kriegerfamily} follow  because
$\{t_{\lambda}:\lambda\in\overline{\Lambda}\}$ is a Cuntz-Krieger
$\overline{\Lambda}$-family.  Lemma \ref{lem:lambdamin} implies that
$\Lbarmin(\lambda,\mu)\subseteq \Lambda$, which shows Condition
(TCK3) is satisfied.  Lemma \ref{lem:finiteexhaustive} gives that
any finite exhaustive subset of $\Lambda$ is a finite exhaustive
subset of $\overline{\Lambda}$.  Therefore, the fact that
$\{t_{\lambda}:\lambda\in\overline{\Lambda}\}$ is a Cuntz-Krieger
$\overline{\Lambda}$-family implies that Condition (CK) of Definition
\ref{def:Toeplitz-Cuntz-Kriegerfamily} is satisfied, proving the
result.
\end{proof}

In the next theorem, we show that $\cst(\Lambda)$ is naturally isomorphic to a subalgebra of $\cst(\overline{\Lambda})$.  The isomorphism is natural in the sense that $\cst(\Lambda)$ is isomorphic to the $\cst$-algebra generated by the set of elements of the form $t_\lambda$ where $\lambda$ is a path in the original $k$-graph, $\Lambda$.  Furthermore, the isomorphism maps generators to elements in the canonical way.

\begin{theorem}\label{lambdaiso}
Let $(\Lambda,d)$ be a finitely aligned $k$-graph and let $(\overline{\Lambda},\overline{d})$ be the $k$-graph given in Definition~\ref{def:Lambdabar}.  Let
$\{t_\lambda:\lambda\in\overline{\Lambda}\}$ be a Cuntz-Krieger
$\overline{\Lambda}$ family.  Then $C^*(\Lambda)$ is isomorphic to the subalgebra of $\cst(\overline{\Lambda})$ generated by the set $\{t_\lambda:\lambda\in\Lambda\}$.
\end{theorem}

\begin{proof}
Let $C^*(\overline{\Lambda})$ be  generated by
$\{t_\lambda:\lambda\in\overline{\Lambda}\}$, and let $C^*(\Lambda)$
be generated by the Cuntz-Krieger $\Lambda$-family
$\{s_\lambda:\lambda\in\Lambda\}$.  Let
$A=C^*(\{t_\lambda:\lambda\in\Lambda\})\subseteq
C^*(\overline{\Lambda})$.   By Theorem \ref{containsCKlambdafamily},
$\{t_\lambda:\lambda\in\Lambda\}$ is a Cuntz-Krieger
$\Lambda$-family; thus the universal property of $C^*(\Lambda)$
gives a *-homomorphism $\pi:C^*(\Lambda)\to C^*(\overline{\Lambda})$
such that $\pi(s_\lambda)=t_\lambda$ for all $\lambda\in\Lambda$.
Since
$\{t_\lambda:\lambda\in\Lambda\}=\{\pi(s_\lambda):\lambda\in\Lambda\}$,
it follows that $\pi(C^*(\Lambda))\subseteq A$.  Furthermore,
because $\pi$ maps $C^*(\Lambda)$ onto the set of generators of $A$,
we have $A\subseteq \pi(C^*(\Lambda))$.  Therefore
$\pi(C^*(\Lambda))=A$.  Since $t_v\not=0$ for all
$v\in\Lambda^0\subseteq\overline{\Lambda}^0$, it follows that
$\pi(s_v)=t_v\not=0$ for all $v\in\Lambda^0$.

Let $\theta:\T^k\to \aut(C^*(\overline{\Lambda}))$ denote the gauge
action on $C^*(\overline{\Lambda})$ and $\gamma:\T^k\to
\aut(C^*(\Lambda))$ denote the gauge action on $C^*(\Lambda)$.  For
all $z\in\T^k$ and $\lambda,\mu\in\Lambda$,
\begin{align*}
(\theta_z\circ\pi)(s_\lambda s_\mu^*)&=\theta_z(t_\lambda t_\mu^*)\\
&=z^{d(\lambda)-d(\mu)}t_\lambda t_\mu^*\\
&=\pi(z^{d(\lambda)-d(\mu)}s_\lambda
s_\mu^*)=(\pi\circ\gamma_z)(s_\lambda s_\mu^*).
\end{align*}
It follows then that $\theta_z\circ\pi=\pi\circ\gamma_z$ for all
$z\in\T^k$.   Therefore by \cite[Theorem~4.2]{RSY2}, $\pi$
is injective.  The previous paragraph shows that $\pi$ maps
$C^*(\Lambda)$ surjectively onto $A$. Thus $C^*(\Lambda)\cong A$.
\end{proof}

\begin{theorem}\label{fullcorner}
Let $(\Lambda,d)$ be a row-finite $k$-graph and let $(\overline{\Lambda},\overline{d})$ be the $k$-graph given in Definition~\ref{def:Lambdabar}.   Then $C^*(\Lambda)$ is
a full corner of $C^*(\overline{\Lambda})$.
\end{theorem}

It is in the following proof that the row-finite condition of $\Lambda$ is necessary.  The row-finiteness of $\Lambda$ implies that its extension, $\overline{\Lambda}$ is also row-finite and does not have any sources.  Thus, there are two equivalent sets of Cuntz-Krieger relations that can be used to define $\cst(\overline{\Lambda})$.  In the proof of Theorem
\ref{fullcorner} we use both Condition (CK) of Definition
\ref{def:Toeplitz-Cuntz-Kriegerfamily} and Condition (CK'), which is stated in
Remark~\ref{rmk:locallyconvexCKrelation}.

\begin{proof}
Suppose $\cst(\overline{\Lambda})$ is generated by
$\{t_\lambda:\lambda\in\overline{\Lambda}\}$.  Let
$A=C^*(\{t_\lambda:\ \lambda\in\Lambda\})\subseteq
\cst(\overline{\Lambda})$. Then $A\cong C^*(\Lambda)$ by Theorem
\ref{lambdaiso}.  We will show that $A$ is a full corner of
$C^*(\overline{\Lambda})$.

Using an argument like that in \cite[Lemma 1.29(c)]{BPRS},
$\sum_{v\in\Lambda^0} t_v$ converges  strictly in
$M(C^*(\overline{\Lambda}))$ to a projection $p$ satisfying
$$p t_\lambda t_{\mu}^* p = \begin{cases}
t_\lambda t_{\mu}^* &\text{if } \overline{r}(\lambda),\overline{r}(\mu)\in \Lambda^0,\\
0 &\text{otherwise.}
\end{cases}$$

Therefore, for  all $\lambda,\mu \in \Lambda$, $t_\lambda
t_{\mu}^*=p t_\lambda t_{\mu}^* p \in pC^*(\overline{\Lambda})p$.
Hence $A\subseteq pC^*(\overline{\Lambda})p$.

If either $\overline{r}(\lambda)$ or $\overline{r}(\mu)$ is  in
$\overline{\Lambda}^0\backslash \Lambda$, then $p t_\lambda
t_{\mu}^* p=0$.  Furthermore, if
$\overline{s}(\lambda)\not=\overline{s}(\mu)$, then $t_\lambda
t_{\mu}^*=0=p t_\lambda t_{\mu}^* p$.  Suppose
$\lambda,\mu\in\overline{\Lambda}$ such that
$\overline{r}(\lambda),\overline{r}(\mu)\in\Lambda^0$ and
$\overline{s}(\lambda)=\overline{s}(\mu)$.
\\

\textsc{Claim:} If $\lambda,\mu\in\overline{\Lambda}$  with
$\overline{r}(\lambda),\overline{r}(\mu)\in\Lambda^0$ and
$\overline{s}(\lambda)=\overline{s}(\mu)\not\in\Lambda^0$, then
$pt_{\lambda}t_{\mu}^* p$ is an element of $A$.
\\

\emph{Proof of Claim:} There exist $x,y\in\Lambda^{\leq\infty},$ and
$l,m,n,q\in\N^k$  such that $\lambda=\pec(x,l,m)$ and
$\mu=\pec(y,n,q)$.  Without loss of generality we may assume
$l=n=0$.  We will proceed by induction on $m$.

Notice that if $m=m\wedge d(x)$, then $\lambda,\mu\in\Lambda$.
Thus~$pt_{\lambda}t_{\mu}^*p=t_{\lambda}t_{\mu}^*\in A$~.

Suppose $m> m\wedge d(x)$, and suppose for an inductive  hypothesis
that the Claim holds for all $n<m$ such that
$\overline{s}(\lambda)=\overline{s}(\mu)=\vcc(x,n)$.

Since $\overline{s}(\lambda)=\overline{s}(\mu)$, (V1) and (V2) of
Definition \ref{def:vertexequivalence} imply that $x(m\wedge
d(x))=y(q\wedge d(y))$ and $m-m\wedge d(x)=q-q\wedge d(y)$.
Therefore $\pec(x,{m\wedge d(x)},m)=\pec(y,{q\wedge d(y)},q)$.  Let
\begin{gather*}
\lambda'=x(0,m\wedge d(x)),\\
\mu'=y(0,q\wedge d(y)), \text{ and}\\
\nu=\pec(x,{m\wedge d(x)},m)=\pec(y,{q\wedge d(y)},q).
\end{gather*}
Then $\lambda=\lambda'\nu$ and $\mu=\mu'\nu$.  There are two cases
to consider.
\\

\textsc{Case 1:} There exist $i_0,i_1\in \{1,2,\ldots,k\}$ such that
$m_{i_j}\geq d(x)_{i_j}+1$.

Let $a=m-e_{i_0}$.  Then $m\wedge d(x)<a<m$, and $a\wedge
d(x)=m\wedge d(x)$. Furthermore,  $\nu=\pec(x,{m\wedge
d(x)},a)\pec(x,a,m)$.  We claim that $\{\pec(x,a,m)\}$ is a finite
exhaustive subset of $\vcc(x,a)\overline{\Lambda}$.  Suppose
$\pec(z,t,u)\in \vcc(x,a)\overline{\Lambda}$.  Then
$\pec(z,t,{t+(m-a)\vee (u-t)})$ is a minimal common extension of
$\pec(z,t,u)$ and $\pec(x,a,m)$.  To see this, we must show that
$\pec(z,t,t+m-a)=\pec(x,a,m)$.  Since
$\vcc(z,t)=\overline{r}(\pec(z,t,u))=\vcc(x,a)$ it follows that
\begin{equation}\label{fullcornereq1}
 z(t\wedge d(z))=x(a\wedge d(x)) \text { and } t-t\wedge d(z)=a-a\wedge d(x).
\end{equation}

Since $a_i=m_i$ for $i\not=i_0$, we have $t_i+m_i-a_i=t_i$ and  so
$((t+m-a)\wedge d(z))_i=(t\wedge d(z))_i$ if $i\not=i_0$.  Since
$m_{i_0}-a_{i_0}=1$ and $m_{i_0}\geq d(x)_{i_0}+1$, it follows that
$a_{i_0}\geq d(x)_{i_0}$ which implies that $t_{i_0}\geq d(z)_{i_0}$
because of \eqref{fullcornereq1}.  Thus $d(z)_{i_0}\leq
t_{i_0}<t_{i_0}+1=t_{i_0}+m_{i_0}-a_{i_0}$ which implies
$d(z)_{i_0}=(t\wedge d(z))_{i_0}=((t+m-a)\wedge d(z))_{i_0}$.  Hence
\begin{align*}
z(t\wedge d(z),(t+m-a)\wedge d(z))&=z(t\wedge d(z),t\wedge d(z))\\
&=z(t\wedge d(z))\\
&=x(a\wedge d(x))\\
&=x(a\wedge d(x),m\wedge d(x))
\end{align*}
because $a\wedge d(x)=m\wedge d(x)$. By \eqref{fullcornereq1} and
the  fact that $t+m-a-t=m-a$, it follows that
$\pec(z,t,t+m-a)=\pec(x,a,m)$.  Therefore, we obtain  that
$$(\pec(z,u,{t+(m-a)\vee(u-t)}),\pec(z,{t+m-a},{t+(m-a)\vee (u-t)}))\in \overline{\Lambda}^{min}(\pec(z,t,u),\pec(x,a,m)).$$  Since $\pec(z,t,u)\in \vcc(x,a)\overline{\Lambda}$ was arbitrary, this implies that $\{\pec(x,a,m)\}$ is a finite exhaustive subset of $\vcc(x,a)\overline{\Lambda}$.

Let $\nu'=\pec(x,m\wedge d(x),a)$.  Then $\nu=\nu'\pec(x,a,m)$, and $\overline{r}(\nu')=\overline{r}(\nu)=\overline{s}(\lambda')$.  Furthermore
\begin{align*}
p t_{\lambda}t_{\mu}^* p&= p t_{\lambda'\nu}t_{\mu'\nu}^* p\\
&= p t_{\lambda'}t_{\nu}t_{\nu}^*t_{\mu'}^* p \\
&= p t_{\lambda'}t_{\nu'\pec(x,a,m)}t_{\nu'\pec(x,a,m)}^*t_{\mu'}^* p\\
&= p t_{\lambda'}t_{\nu'}t_{\pec(x,a,m)}t_{\pec(x,a,m)}^*t_{\nu'}^*t_{\mu'}^* p\\
&= p t_{\lambda'}t_{\nu'}t_{\vcc(x,a)}t_{\nu'}^*t_{\mu'}^* p \text{ because $\{\pec(x,a,m)\}\in \vcc(x,a)\FE(\overline{\Lambda})$}\\
&= p t_{\lambda'}t_{\nu'}t_{\nu'}^*t_{\mu'}^* p\\
&= p t_{\lambda'\nu'}t_{\mu'\nu'}^* p\\
\end{align*}
which belongs to $A$ by the inductive hypothesis since
$\overline{s}(\nu')=\vcc(x,a)$ and $a<m$.   This concludes Case 1.
\\

\textsc{Case 2:} Suppose that $m=m\wedge d(x)+e_{i_0}$ for some
$i_0\in \{1,2,\ldots,k\}$.   Let $v$ be the vertex $x(m\wedge d(x))$.  We will show that $v\overline{\Lambda}^{e_{i_0}}\backslash\Lambda$ is the set $\{\nu\}$.
Let $\xi\not\in\Lambda$ be an element of $v\overline{\Lambda}^{e_{i_0}}$. Then $\xi=\pec(z,0,e_{i_0})$ for some $z\in
v\Lambda^{\leq\infty}$. Since $e_{i_0}\not\leq d(z)$ and
$(e_{i_0})_j\leq d(z)_j$ for $j\not=i_0$, it must be that
$d(z)_{i_0}=0$.  Then $e_{i_0}\wedge d(z)=0$; $0-(e_{i_0}\wedge
d(z)=0$ and $z(0,e_{i_0}\wedge d(z))=z(0,0)=v$.  Since $m-m\wedge
d(x)=e_{i_0}$ and $0=m\wedge d(x)-(m\wedge d(x))\wedge d(x)$, it follows that  $\pec(z,0,e_{i_0})=\nu$, which is a contradiction.

Let $E=v\overline{\Lambda}^{e_{i_0}}\bigcap \Lambda$.  Then
$E=v\overline{\Lambda}^{e_{i_0}}\backslash\{\nu\}$.  Since
$\overline{\Lambda}$ has no
sources by Theorem \ref{Lambdabarproperties} and is row-finite by Theorem~\ref{thm:finitelyalignedpreservedinLbar}, we have that
$v\overline{\Lambda}^{\leq e_{i_0}}=v\overline{\Lambda}^{e_{i_0}}$.
Then by \cite[Proposition B.1]{RSY2},
$$t_v=\sum_{\xi\in v\overline{\Lambda}^{\leq e_{i_0}}} t_\xi t_\xi^*= t_{\nu}t_{\nu}^*+\sum_{\lambda\in E} t_\xi t_\xi^* .$$
Thus
\begin{align*}
p t_{\lambda}t_{\mu}^* p&= p t_{\lambda'\nu}t_{\mu'\nu}^* p\\
&= p t_{\lambda'}t_{\nu}t_{\nu}^*t_{\mu'}^* p \\
&= p t_{\lambda'}(t_v-\sum_{\xi\in E} t_\xi t_\xi^*)t_{\mu'}^* p\\
\end{align*}
which belongs to $A$ because $v\in\Lambda$ and $E\subseteq\Lambda$.
This concludes Case 2,  and proves the claim.
\\

Therefore $p C^*(\overline{\Lambda})p\subseteq A$.  Hence $A=p
C^*(\overline{\Lambda})p$.

To show that $A$ is a full corner of $C^*(\overline{\Lambda})$,
suppose that $J$ is an ideal  in $C^*(\overline{\Lambda})$ such that
$A\subseteq J$.  Of course $\{t_\lambda:\lambda\in\Lambda\}\subseteq
J$ because this set generates $A$.  Let
$v\in\overline{\Lambda}^0\backslash\Lambda$. Then $v=\vcc(x,m)$ for
some $x\in\Lambda^{\leq\infty}$ and $m\not\leq d(x)$.  Then the path
$\alpha=\pec(x,m\wedge d(x),m)\in\overline{\Lambda}$ and
$\overline{r}(\alpha)=x(m\wedge d(x))\in\Lambda^0$.  Also
$\overline{s}(\alpha)=\vcc(x,m)$.  Thus $t_\alpha=t_{x(m\wedge
d(x))} t_\alpha \in J$ because $t_{x(m\wedge d(x))}\in J$. Therefore
$t_{\vcc(x,m)}=t_\alpha^* t_\alpha\in J$ and $\{t_v:
v\in\overline{\Lambda}^0\}\subseteq J$.  Next let
$\lambda\in\overline{\Lambda}\backslash\overline{\Lambda}^0$.  Then
$\overline{r}(\lambda)\in\overline{\Lambda}^0$ and
$t_\lambda=t_{\overline{r}(\lambda)}t_\lambda\in J$.  Hence
$\{t_\lambda:\lambda \in\overline{\Lambda}\}$, the set of generators
of $C^*(\overline{\Lambda})$ lies in $J$, which implies that
$J=C^*(\overline{\Lambda})$.
\end{proof}

We now conclude the chapter with the proof of Theorem
\ref{thm:nosources}.

\begin{proof}[Proof of Theorem \ref{thm:nosources}:]
The pair $(\overline{\Lambda},\overline{d})$ of Definition~\ref{def:Lambdabar} is a row-finite
$k$-graph without sources by  Theorems~\ref{thm:finitelyalignedpreservedinLbar} and \ref{Lambdabarproperties}.  By
definition of $\overline{\Lambda}$,
$\Obj(\Lambda)\subseteq\Obj(\overline{\Lambda})$, and
$\Hom(\Lambda)\subseteq\Hom(\overline{\Lambda})$. Furthermore,
$\overline{r}|_{\Hom(\Lambda)}=r$,
$\overline{s}|_{\Hom(\Lambda)}=s$, and $\overline{d}|_{\Lambda}=d$.
Thus the map $\iota:\Lambda\to\overline{\Lambda}$, given by
$\iota(\lambda)=\lambda$ for all $\lambda\in\Lambda$ is a $k$-graph
isomorphism between $\Lambda$ and $\iota\Lambda$.  Therefore
$A=\{t_\lambda:\lambda\in\iota\Lambda\}$ is isomorphic to
$C^*(\Lambda)$ by Theorem \ref{lambdaiso}, and is a full corner of
$C^*(\overline{\Lambda})$ by Theorem \ref{fullcorner}.
\end{proof}

\section{Examples}\label{sec:examples}
In this section, we  will apply the construction of Section
\ref{sec:removing sources} to several examples of row-finite
$k$-graphs. The examples include $k$-graphs that are and are not
locally convex. The examples were chosen to illustrate how the
conditions in Definitions \ref{def:vertexequivalence} and
\ref{def:pathequivalence} affect the construction as well as why
they are necessary.  For the diagrams in this chapter, edges of
degree $(1,0)$ appearing in the original $k$-graph will be drawn
with double solid arrows (\xy\ar@{=>}(0,0);(5,0)\endxy); edges of
degree $(0,1)$ in the original $k$-graph will be drawn with double
dashed arrows (\xy\ar@{==>} (0,0);(5,0)\endxy).  Edges of degree
$(1,0)$ and $(0,1)$ that appear in the extension will be
represented, respectively, by solid arrows (\xy\ar(0,0);(5,0)\endxy)
and dashed arrows (\xy\ar@{-->}(0,0);(5,0)\endxy).

\begin{example}\label{ex 1}
Let $\Lambda$  be a row-finite 1-graph with sources.  In \cite{BPRS}
and \cite{DT} the method of ``adding heads to sources" was used to
create a row-finite 1-graph without sources that preserved the
Morita equivalence class of $C^*(\Lambda)$.  We will show that the
method developed in Chapter 3 coincides with the previous
construction of \cite{BPRS,DT}.

Let $\Lambda_S=\{v\in\Lambda^0: v\Lambda^1=\emptyset\}$.   Let
$v\in\Lambda_S$.  Then $\Lambda_S$ is the set of sources as defined
for a directed graph.  In \cite{BPRS}, \emph{adding a head to $v$}
means attaching the following graph  to $v$.
$$
\xymatrix{v&v_1\ar[l]_{e_{v_1}}&v_2\ar[l]_{e_{v_2}}&v_3\ar[l]_{e_{v_3}}&\cdots\ar[l]&v_{n-1}\ar[l]&v_{n}\ar[l]_{e_{v_n}}&\cdots\ar[l]}
$$

Let $\Gamma$ denote the 1-graph that results from adding a head to
each $v\in\Lambda_S$. Then any path in $\Gamma$ is either a path in
$\Lambda$ or it is of the form $\lambda e_{v_1}e_{v_2}\ldots
e_{v_n}$ for some $v\in\Lambda_S$, $\lambda\in \Lambda v$ and
$n\in\N$ with $n\geq 1$.

Suppose $x:\Omega_{1,m}\to\Lambda$ is a graph  morphism for some
$m\in \N$ (so we are considering only finite paths).   Then
$x\in\Linf$ if and only if  $x(m)=x(d(x))\in\Lambda_S$. Thus,
\begin{gather*}
V_\Lambda=\bigcup_{v\in\Lambda_S}\{(x;m):x(d(x))=v \text{ and }
m>d(x)\}  \text{ and}\\
P_\Lambda=\bigcup_{v\in\Lambda_S}\{(x;(m,n)):x(d(x))=v, m\leq n
\text{ and } n>d(x)\}.
\end{gather*}

Suppose $x$ and $y$ are paths in $\Linf$ such that $d(x)$ and $d(y)$
are finite.   Suppose further that $x(d(x))=y(d(y))=v$ for some
source $v\in\Lambda_S$.  Since
$\shift({d(x)},x)=\shift({d(y)},y)=v$, Proposition
\ref{prop:shiftequivalence} implies
$\vcc(x,{d(x)+m})=\vcc(y,{d(y)+m})=\vcc(v,m)$ for all $m\in\N$,
$m\geq 1$.  Also by Proposition \ref{prop:shiftequivalence}, for
$m,n\in\N$ with $m\leq n$, we have
$$\pec(x,{d(x)+m},{d(x)+n})=\pec(y,{d(y)+m},{d(y)+n})=\pec(v,m,n).$$   So  for any
$\pec(x,m,n)\in\widetilde{P_\Lambda}$, let $v_x=x(d(x))$.  Then
$v_x\in\Lambda_S$ and we have that
$$\pec(x,m,n)=\begin{cases}
\vcc(v_x,{m-d(x)},{n-d(x)})&\quad\text{if } m\geq d(x),\\
x(m,d(x))\pec(v_x,0,{n-d(x)})&\quad\text{if } m< d(x).
\end{cases}$$
Therefore, the vertices and paths added to $\Lambda$ to form $\Lbar$
are
\begin{gather*}
\widetilde{V_\Lambda}=\bigcup_{v\in\Lambda_S}\{\vcc(v,m):m\geq 1\}, \text{ and}\\
\widetilde{P_\Lambda}=\bigcup_{v\in\Lambda_S}\{\pec(v,m,n):m,n\in\N,
m\leq n\}\cup\{\lambda\pec(v,0,n):\lambda\in \Lambda v,n>0\}.
\end{gather*}

The assignment $\vcc(v,m)\mapsto v_m$ and $\pec(v,m-1,m)\mapsto
e_{v_{m}}$ for  all $v\in\Lambda^0$ and $m\in\N$ with $m\geq 1$
creates a graph isomorphism between $\overline{\Lambda}$ and
$\Gamma$ when it is extended in a natural way to the entire
category.  That is, define $\Phi:\overline{\Lambda}\to \Gamma$ by
\begin{align*}
&\Phi(\lambda)=\lambda \text{ for all } \lambda\in\Lambda,\\
&\Phi(\vcc(v,m))=v_m \text{ for all } \vcc(v,m)\in\Obj(\overline{\Lambda}),\\
&\Phi(\pec(v,m,n))=e_{v_{m+1}}e_{v_{m+2}}\ldots e_{v_n} \text{ for all } v\in\Lambda_S, m\leq n,\text{ and }\\
&\Phi(\lambda\pec(v,0,n))=\lambda e_{v_1}e_{v_2}\ldots e_{v_n}
\text{ for all } v\in\Lambda_S, \lambda\in v\Lambda, n\in\N.
\end{align*}
Then $\Phi$ is a graph isomorphism, and so for 1-graphs, the
desingularization  developed in Section \ref{sec:removing sources}
is the same as the method used in \cite{BPRS, DT}.
\end{example}

\begin{example}$\mathbf{\Omega_{k,m}.}$\label{ex2}
Let $\Lambda$ be the 2-graph $\Omega_{2,(1,1)}$ shown below.
$$
\xymatrix{v_2\ar@{==>}[d]_{\mu}&v_3\ar@{=>}[l]_{\beta}\ar@{==>}[d]^{\alpha}\\
v_0&v_1\ar@{=>}[l]^{\lambda}}
$$

For this example, $\Linf$ consists of four elements:
$$\begin{array}{lll}
w:\Omega_{2,(0,0)}\to\Lambda&\quad&x:\Omega_{2,(0,1)}\to\Lambda\\
w((0,0))=v_3&&x((0,0),(0,1))=\alpha\\
\\
y:\Omega_{2,(1,0)}\to\Lambda&\quad&z:\Omega_{2,(1,1)}\to\Lambda\\
y((0,0),(1,0))=\beta&&z((0,0),(1,1))=\lambda\alpha=\mu\beta
\end{array}$$

Since $d(w)=(0,0)$, the set $\{\vcc(w,m):m\in\N^2, m>(0,0)\}$ lies  in
$\widetilde{V_\Lambda}$, and the set $\{\pec(w,m,n): m,n\in\N^2
\text{ and }m\leq n\}$ is a subset of $\widetilde{X_{\Lambda}}$. The
figure that follows shows the 1-skeleton of these elements together
with the original graph $\Lambda$.  In this figure,
$a_w=\vcc(w,{(1,1)})$, $b_w=\vcc(w,{(1,2)})$ and
$\xi_w=\pec(w,{(1,1)},{(1,2)})$.

$$
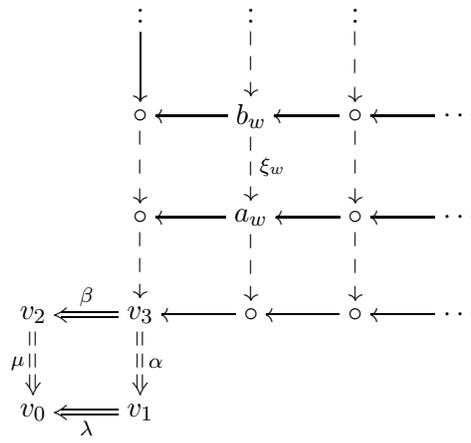
\xymatrix{&\vdots\ar@{->}[d]&\vdots\ar@{-->}[d]&\vdots\ar@{-->}[d]\\
&\circ\ar@{-->}[d]&b_w\ar@{-->}[d]^{\xi_w}\ar[l]&\circ\ar@{-->}[d]\ar[l]&\cdots\ar[l]\\
&\circ\ar@{-->}[d]&a_w\ar@{-->}[d]\ar[l]&\circ\ar@{-->}[d]\ar[l]&\cdots\ar[l]\\
v_2\ar@{==>}[d]_{\mu}&v_3\ar@{=>}[l]_{\beta}\ar@{==>}[d]^{\alpha}&\circ\ar[l]&\circ\ar[l]&\cdots\ar[l]\\
v_0&v_1\ar@{=>}[l]^{\lambda}}$$

From the boundary path $x$,  we have $\{\vcc(x,m):m\in\N^2,
m\not\leq (0,1)\}\subset \widetilde{V_\Lambda}$ and $\{\pec(x,m,n)):
m\leq n, n\not\leq (0,1)\}\subset \widetilde{P_\Lambda}$.  Below, we
see the 1-skeleton of these elements as well as $\Lambda$. Here
$a_x=\vcc(x,{(1,2)})$, $b_x=\vcc(x,{(1,3)})$ and
$\xi_x=\pec(x,{(1,2)},{(1,3)})$.

$$
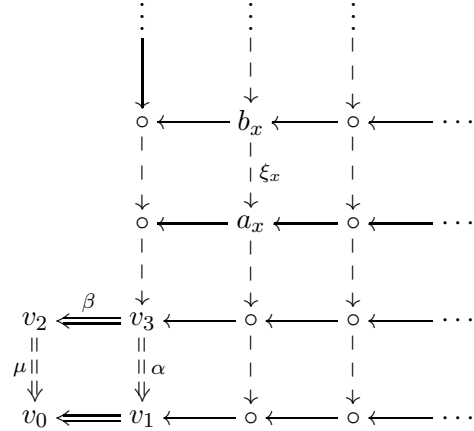
\xymatrix{&\vdots\ar@{->}[d]&\vdots\ar@{-->}[d]&\vdots\ar@{-->}[d]\\
&\circ\ar@{-->}[d]&b_x\ar@{-->}[d]^{\xi_x}\ar[l]&\circ\ar@{-->}[d]\ar[l]&\cdots\ar[l]\\
&\circ\ar@{-->}[d]&a_x\ar@{-->}[d]\ar[l]&\circ\ar@{-->}[d]\ar[l]&\cdots\ar[l]\\
v_2\ar@{==>}[d]_{\mu}&v_3\ar@{=>}[l]_{\beta}\ar@{==>}[d]^{\alpha}&\circ\ar@{-->}[d]\ar[l]&\circ\ar@{-->}[d]\ar[l]&\cdots\ar[l]\\
v_0&v_1\ar@{=>}[l]&\circ\ar[l]&\circ\ar[l]&\cdots\ar[l]}$$

The elements of $V_\Lambda$ and $P_\Lambda$ resulting from the
boundary paths $y$ and $z$ are similar.  The next two figure show
$\Lambda$ together with the additional vertices and paths. In the
first figure that follows, we have $a_y=\vcc(y,{(2,1)})$,
$b_y=\vcc(y,{(2,2)})$ and $\xi_y=\pec(y,{(2,1)},{(2,2)})$, while
$a_z=\vcc(z,{(2,2)})$, $b_z=\vcc(z,{(2,3)})$ and
$\xi_z=\pec(z,{(2,2)},{(2,3)})$ in the second.

$$
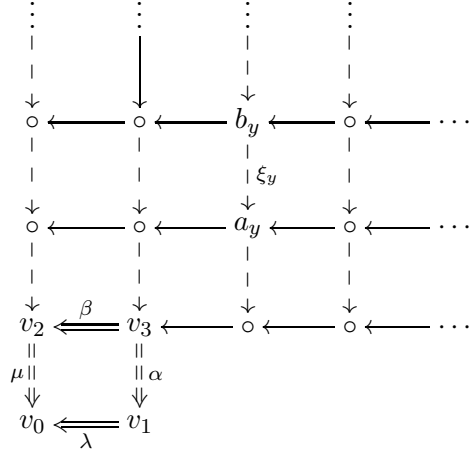
\xymatrix{\vdots\ar@{-->}[d]&\vdots\ar@{->}[d]&\vdots\ar@{-->}[d]&\vdots\ar@{-->}[d]\\
\circ\ar@{-->}[d]&\circ\ar@{-->}[d]\ar[l]&b_y\ar@{-->}[d]^{\xi_y}\ar[l]&\circ\ar@{-->}[d]\ar[l]&\cdots\ar[l]\\
\circ\ar@{-->}[d]&\circ\ar@{-->}[d]\ar[l]&a_y\ar@{-->}[d]\ar[l]&\circ\ar@{-->}[d]\ar[l]&\cdots\ar[l]\\
v_2\ar@{==>}[d]_{\mu}&v_3\ar@{=>}[l]_{\beta}\ar@{==>}[d]^{\alpha}&\circ\ar[l]&\circ\ar[l]&\cdots\ar[l]\\
v_0&v_1\ar@{=>}[l]^{\lambda}}
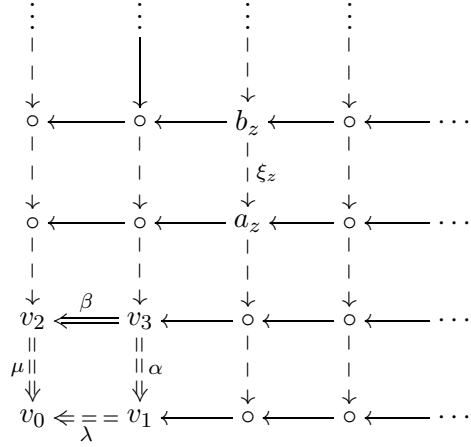
$$

$$
\xymatrix{\vdots\ar@{-->}[d]&\vdots\ar@{->}[d]&\vdots\ar@{-->}[d]&\vdots\ar@{-->}[d]\\
\circ\ar@{-->}[d]&\circ\ar@{-->}[d]\ar[l]&b_z\ar@{-->}[d]^{\xi_z}\ar[l]&\circ\ar@{-->}[d]\ar[l]&\cdots\ar[l]\\
\circ\ar@{-->}[d]&\circ\ar@{-->}[d]\ar[l]&a_z\ar@{-->}[d]\ar[l]&\circ\ar@{-->}[d]\ar[l]&\cdots\ar[l]\\
v_2\ar@{==>}[d]_{\mu}&v_3\ar@{=>}[l]_{\beta}\ar@{==>}[d]^{\alpha}&\circ\ar@{-->}[d]\ar[l]&\circ\ar@{-->}[d]\ar[l]&\cdots\ar[l]\\
v_0&v_1\ar@{==>}[l]^{\lambda}&\circ\ar[l]&\circ\ar[l]&\cdots\ar[l]}$$

Since $x=\shift({(1,0)},z)$, Proposition \ref{prop:shiftequivalence}
implies that $\vcc(x,m)=\vcc(z,{m+(1,0)})$ for all $m\not\leq(1,0)$,
and $\pec(x,m,n)=\pec(z,{m+(1,0)},{n+(1,0)})$ for all $m\leq n$,
$n\not\leq (1,0)$.  Similarly $y=\shift({(0,1)},z)$ and
$w=\shift({(1,1)},w)$.  Therefore by
Proposition~\ref{prop:shiftequivalence}, we obtain the following equalities
\begin{align*}
\vcc(y,m)=\vcc(z,{m+(0,1)})& \text{ for all } m\not\leq (0,1);\\
\pec(y,m,n)=\pec(z,{m+(0,1)},{n+(0,1)})& \text{ for all } m\leq n, n\not\leq (0,1);\\
\vcc(w,m)=\vcc(z,{m+(1,1)})& \text{ for all } m>0, \text{ and} \\
\pec(w,m,n)=\pec(z,{m+(1,1)},{n+(1,1)})& \text{ for all } m\leq n,
n>0.
\end{align*}

Thus,
\begin{gather*}
\widetilde{V_\Lambda}=\{\vcc(z,m):m\not\leq (1,1)\},\text{ and}\\
\widetilde{P_\Lambda}=\{\pec(z,m,n):m\leq n \text{ and } n\not\leq
(1,1)\}.
\end{gather*}
Therefore, $\overline{\Lambda}$ is $\Omega_{2,(\infty,\infty)}$.

It can be shown that $\cst(\Lambda)\cong M_4(\C)$ and  that
$\cst(\Lbar)\cong\mathcal{K}(\ell^2(\N^2))$.  So we see that
$\cst(\Lambda)$ is indeed a full corner of $\cst(\Lbar)$.

In general, if $\Lambda=\Omega_{k,m}$ for  some $m\in
(\N\cup\{\infty\})^k$, then $\overline{\Lambda}=\Omega_{k}$. This
seems reasonable since $\Omega_k$ is the simplest $k$-graph without
sources that contains $\Omega_{k,m}$ as a subgraph.  In a sense, we
are just ``filling in the gaps" of $\Omega_{k,m}$ to extend it to
$\Omega_k$.
\end{example}

\begin{example}\textbf{A non-locally convex graph.}
Let $\Lambda$ be the 2-graph shown below.

$$
\xymatrix{v_2\ar@{==>}[d]_{\mu}\\
v_0&v_1\ar@{=>}[l]^{\lambda}}
$$

While $\Lambda$ is a subgraph of $\Omega_{2,(\infty,\infty)}$, the $\cst$-algebra of $\Lambda$ will not sit inside $\cst(\Omega_{2,(\infty,\infty)})$ as a full corner.  According to \cite{S2}, $\cst(\Lambda)$ will have two maximal ideals corresponding to the saturated and hereditary subsets of $\Lambda$ which are $\{v_1\}$ and $\{v_2\}$.  However, $\cst(\Omega_{2,(\infty,\infty)})$ is a simple $\cst$-algebra.

For this example, $\Linf$ consists of four boundary paths, but there are only two boundary paths that we must consider.  All other elements of $\Linf$ are shifts of the paths $x$ and $y$ described below.  As in the previous example,
Proposition~\ref{prop:shiftequivalence} implies that $\overline{\Lambda}$ is determined by these paths.

Define $x:\Omega_{2,(1,0)}\to\Lambda$ and $y:\Omega_{2,(0,1)}\to\Lambda$ to be the following graph morphisms.
$$\begin{array}{lll}
x:\Omega_{2,(1,0)}\to\Lambda&\quad&y:\Omega_{2,(0,1)}\to\Lambda\\
x((0,0), (1,0))=\lambda&&y((0,0),(0,1))=\mu.\\
\end{array}$$

Both $x$ and $y$ extend to form a copy of $\Omega_{2,(\infty,\infty)}$ in $\overline{\Lambda}$. However, the extensions of these paths may be equivalent according to Definition~\ref{def:vertexequivalence} or Definition~\ref{def:pathequivalence}.

Let $\vcc(x,m)$ and $\vcc(y,p)$ be elements of $\overline{\Lambda}^0$.  Suppose that $\vcc(x,m)=\vcc(y,p)$.  Then because $x$ and $y$ agree only at $x((0,0))=y((0,0))=v_0$, we must have that $m\wedge d(x)=p\wedge d(y)=(0,0)$ by Condition~(V1) of Definition~\ref{def:vertexequivalence}.  Therefore $m=(0,m_2)$ and $p=(p_1,0)$ for $m_2,p_1>0$.  But Condition~(V2) would imply that $(0,m_2)=(p_1,0)$, which is impossible.  Hence, $\vcc(x,m)\neq\vcc(y,p)$ for all $m\not\leq d(x)$ and $p\not\leq d(y)$.  Hence the two copies of $\Omega_{2,(\infty,\infty)}$ that these boundary paths contribute to $\overline{\Lambda}$ intersect only at $v_0$.  The extension of $\Lambda$ is drawn below.

$$
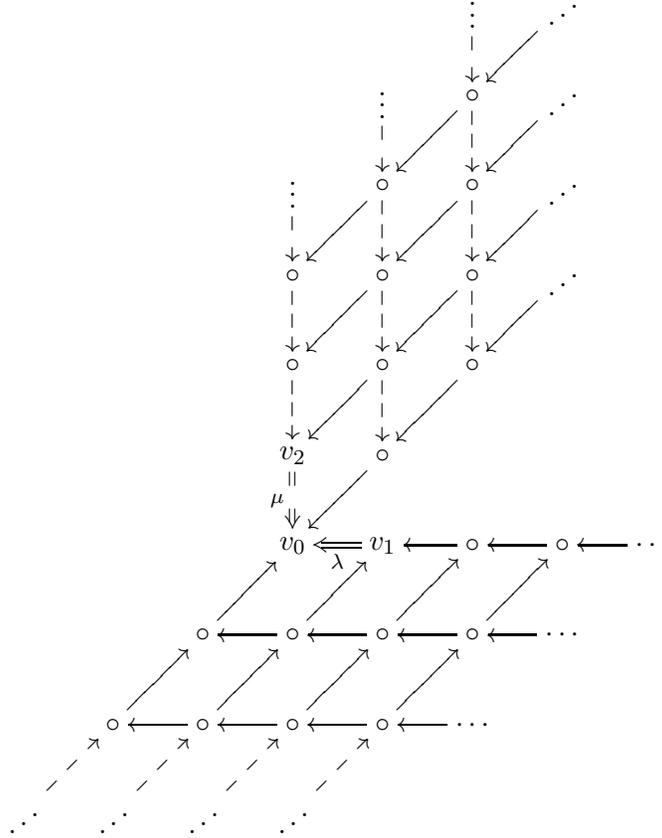
\xymatrix@!R=10pt@!C=10pt{&&&&&\vdots\ar@{-->}[d]&\iddots\ar[dl]\\
&&&&\vdots\ar@{-->}[d]&\circ\ar[dl]\ar@{-->}[d]&\iddots\ar[dl]\\
&&&\vdots\ar@{-->}[d]&\circ\ar[dl]\ar@{-->}[d]&\circ\ar[dl]\ar@{-->}[d]&\iddots\ar[dl]\\
&&&\circ\ar@{-->}[d]&\circ\ar[dl]\ar@{-->}[d]&\circ\ar[dl]\ar@{-->}[d]&\iddots\ar[dl]\\
&&&\circ\ar@{-->}[d]&\circ\ar[dl]\ar@{-->}[d]&\circ\ar[dl]\\
&&&v_2\ar@{==>}[d]_{\mu}&\circ\ar[dl]\\
&&&v_0&v_1\ar@{=>}[l]^{\lambda}&\circ\ar[l]&\circ\ar[l]&\cdots\ar[l]\\
&&\circ\ar@{->}[ur]&\circ\ar@{->}[ur]\ar[l]&\circ\ar@{->}[ur]\ar[l]&\circ\ar@{->}[ur]\ar[l]&\cdots\ar[l]\\
&\circ\ar@{->}[ur]&\circ\ar@{->}[ur]\ar[l]&\circ\ar@{->}[ur]\ar[l]&\circ\ar@{->}[ur]\ar[l]&\cdots\ar[l]\\
\iddots\ar@{-->}[ur]&\iddots\ar@{-->}[ur]&\iddots\ar@{-->}[ur]&\iddots\ar@{-->}[ur]}
$$
\\

For this example $\cst(\Lambda)\cong M_2\oplus M_2$ and $\cst(\overline{\Lambda})\cong \mathcal{K}(\ell^2(\N^2))\oplus \mathcal{K}(\ell^2(\N^2))$.
\end{example}

\section{Additional questions}
For directed graphs, the desingularization process developed in \cite{DT} takes any directed graph with sources and infinite receivers and builds a directed graph without these singular vertices while still preserving the Morita equivalence class of the graph $\cst$-algebras.

Consider the following directed graph $E$.  This graph does not have any sources, but $v$ receives infinitely many edges.  Label the edges from $w$ to $v$ as $\alpha_i$, $i\in\N$.

$$\xymatrix{
\cdots&\bullet\ar[l]&\bullet\ar[l]&v\ar[l]&w\ar[l]_{\infty}&\bullet\ar[l]&\bullet\ar[l]&\cdots\ar[l]}$$

The desingularization process will add a head to $v$ and resdistribute the edges to the new vertices. Let $F$ denote the desingularization of $E$.  The directed graph $F$ is drawn below.  There is a bijection between the set of all finite paths of $E$ and the set of finite paths in $F$ that have range and source in $E$.  This bijection maps $\alpha_1$ to $f_1$ and sends $\alpha_i$, $i>1$ to the path $e_{v_1}e_{v_2}\ldots e_{v_i-1}f_i$.

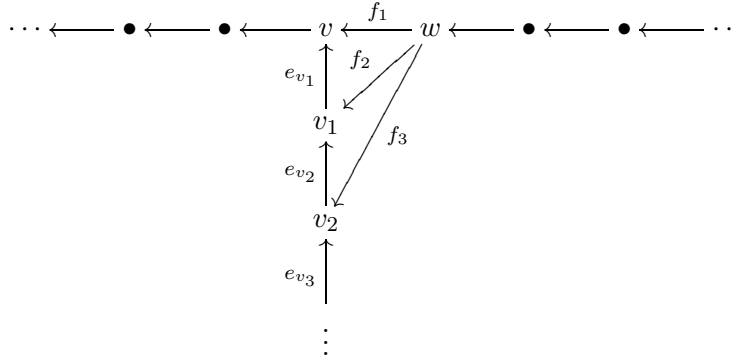
$$\xymatrix{
\cdots&\bullet\ar[l]&\bullet\ar[l]&v\ar[l]&w\ar[l]_{f_1}\ar[dl]_{f_2}\ar[ddl]^{f_3}&\bullet\ar[l]&\bullet\ar[l]&\cdots\ar[l]\\
&&&v_1\ar[u]^{e_{v_1}}\\
&&&v_2\ar[u]^{e_{v_2}}\\
&&&\vdots\ar[u]^{e_{v_3}}\\
}$$

It remains to be seen if a desingularization process for infinite receivers in a higher-rank graph can be developed.  The process outlined in this paper for dealing with sources in a higher-rank graph is analogous to the process of ``adding a head to a source."  When a head is attached to a source in a 1-graph, a copy of $\Omega_{1,\infty}$ is created in the 1-graph.  The method developed in Section~\ref{sec:removing sources}  extends a $k$-graph with sources in a way that creates a copy of $\Omega_{k,(\infty,\ldots,\infty)}$  in the extension.  If the desingularization of a $k$-graph with infinite receivers is to remain analogous to what occurs in the 1-graph setting, then we must redistribute infinitely many edges of various degrees throughout a copy of $\Omega_{k,(\infty,\ldots,\infty)}$.  Deciding how to do this is complicated by the fact that adding just one edge to a vertex often necessitates adding many edges to other vertices to ensure that the factorization property holds.  Furthermore, there are many different ways that a vertex in a $k$-graph can receive infinitely many paths of a certain degree.  For example, in the 2-graphs  $\Lambda_1$ through $\Lambda_4$ below, the vertex $v_0$ receives infinitely many edges of degree $(1,1)$.

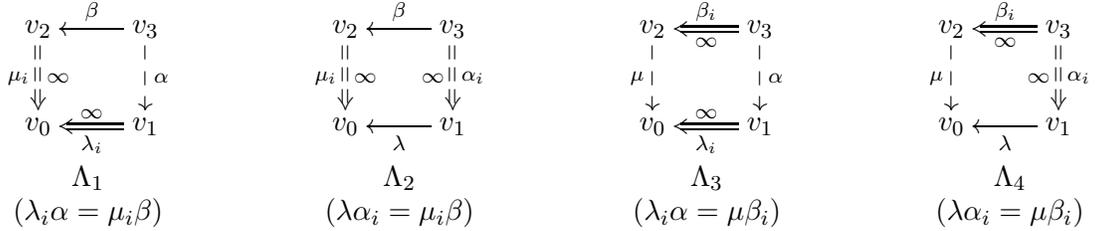
\begin{center}
\begin{tabular}{ccccccc}
$\xymatrix{v_2\ar@{==>}[d]_{\mu_i}^\infty&v_3\ar[l]_{\beta}\ar@{-->}[d]^{\alpha}\\
v_0&v_1\ar@{=>}[l]^{\lambda_i}_\infty}$&\hspace{.4in}\ &
$\xymatrix{v_2\ar@{==>}[d]_{\mu_i}^\infty&v_3\ar[l]_{\beta}\ar@{==>}[d]^{\alpha_i}_\infty\\
v_0&v_1\ar@{->}[l]^{\lambda}}$&\hspace{.4in}\ &
$\xymatrix{v_2\ar@{-->}[d]_{\mu}&v_3\ar@{=>}[l]_{\beta_i}^\infty\ar@{-->}[d]^{\alpha}\\
v_0&v_1\ar@{=>}[l]^{\lambda_i}_\infty}$&\hspace{.4in}\ &
$\xymatrix{v_2\ar@{-->}[d]_{\mu}&v_3\ar@{=>}[l]_{\beta_i}^\infty\ar@{==>}[d]^{\alpha_i}_\infty\\
v_0&v_1\ar@{->}[l]^{\lambda}}$\\
$\Lambda_1$&&$\Lambda_2$&&$\Lambda_3$&&$\Lambda_4$\\
($\lambda_i\alpha=\mu_i\beta$)&&($\lambda\alpha_i=\mu_i\beta$)&&
($\lambda_i\alpha=\mu\beta_i$)&&($\lambda\alpha_i=\mu\beta_i$)
\end{tabular}
\end{center}

When $\Lambda$ is a finitely aligned $k$-graph, the set $\Linf$ is used to create a non-degenerate Toeplitz-Cuntz-Krieger $\Lambda$-family in \cite{RSY1}.  For locally convex, row-finite $k$-graphs, these paths are related to the sets $\Lambda^{\leq n}$, which appear in the Cuntz-Krieger relation (CK$^\prime$) (Remark~\ref{rmk:locallyconvexCKrelation}).  The elements in $\Linf$, in a way, point out where the sources are in the $k$-graph and are crucial to the process developed in Section~\ref{sec:removing sources}.  In \cite{S1}, a different set of boundary paths, the set $\partial\Lambda$, is introduced to study relative Cuntz-Krieger algebras of finitely aligned $k$-graphs.  A graph morphism $x:\Omega_{k,m}\to\Lambda$ belongs to $\partial\Lambda$ if for every $n\leq m$ and every finite-exhaustive set $E\subseteq x(n)\Lambda$, there exists $\mu\in E$ such that $x(n,n+d(\mu))=\mu$ \cite[Definition~4.4]{S1}.  The set $\partial\Lambda$ also plays a part in developing a groupoid model for finitely aligned $k$-graphs \cite{FMY}.  In general, the set $\Linf$ is a proper subset of $\partial\Lambda$, and in some sense, the paths of $\partial\Lambda$ are the limits of sequences of paths in $\Linf$.   The elements in $\partial\Lambda$ identify which vertices in $\Lambda$ are infinite receivers as well as sources.  Perhaps a construction using these paths would lead to a desingularization of a $k$-graph with infinite receivers.


\end{document}